%% file: QMAK_220829.tex
%
\newif\ifloadreferences\loadreferencestrue
%
%
%
%
\let\myfrac=\frac%
\input eplain %
\let\frac=\myfrac%
\let\myfootnote=\footnote%
\input amstex \input epsf %
\let\footnote=\myfootnote%
%
%
\loadeufm\loadmsam\loadmsbm\message{symbol names}\UseAMSsymbols\message{,}%
%
\font\myfontdefault=cmr10%
\newif\ifmakebiblio%
\newif\ifinappendices%
\newif\ifundefinedreferences%
\newif\ifchangedreferences%
\makebibliofalse%
\undefinedreferencesfalse%
\changedreferencesfalse%
%
%
%
%
%
\def\setcatcodes{\catcode`\!=0 \catcode`\\=11}%
{\global\let\noe=\noexpand%
\catcode`\@=11 \catcode`\_=11 \setcatcodes%
!global!def!_@@internal@@makeref#1{%
!global!expandafter!def!csname #1ref!endcsname##1{%
!csname _@#1@##1!endcsname%
!expandafter!ifx!csname _@#1@##1!endcsname!relax%
    !write16{#1 ##1 not defined - run saving references}%
    !undefinedreferencestrue%
!fi}}%
!global!def!_@@internal@@makelabel#1{%
!global!expandafter!def!csname #1label!endcsname##1{%
!edef!temptoken{!csname #1info!endcsname}%
!ifloadreferences%
!expandafter!ifx!csname _@#1@##1!endcsname!relax%
!write16{#1 ##1 not hitherto defined - rerun saving references}%
!changedreferencestrue%
!else%
!expandafter!ifx!csname _@#1@##1!endcsname!temptoken%
!else%
!write16{#1 ##1 reference has changed - rerun saving references}%
!changedreferencestrue%
!fi%
!fi%
!else%
!expandafter!edef!csname _@#1@##1!endcsname{!temptoken}%
!edef!textoutput{!write!references{\global\def\_@#1@##1{!temptoken}}}%
!textoutput%
!fi}}%
!global!def!makecounter#1{!_@@internal@@makelabel{#1}!_@@internal@@makeref{#1}}%
!unsetcatcodes%
}
%
%
%
%
%
\def\turnintolatin#1{\ifcase #1 _\or i\or ii\or iii\or iv\or v\or vi\or vii\or viii\or ix\or x\or xi\or xii\or xiii\or xiv\or xv\or xvi\or xvii\or xviii\or xix\or xx\or xxi\or xxii\or xxiii\or xxiv\or xxv\or xxvi\fi}%
\def\alphanum#1{\ifcase #1 _\or A\or B\or C\or D\or E\or F\or G\or H\or I\or J\or K\or L\or M\or N\or O\or P\or Q\or R\or S\or T\or U\or V\or W\or X\or Y\or Z\fi}%
\newwrite\references%
\ifloadreferences{\catcode`\@=11 \catcode`\_=11 \input references.tex }%
\else{\openout\references=references.tex }%
\fi%
%
%
\newcount\headno%
\global\headno=0%
\def\headinfo{\ifinappendices\alphanum\headno\else\the\headno\fi}%
\def\nextheadno{\global\advance\headno by 1 \global\subheadno=0 \global\eqnno=0 \headinfo}%
\makecounter{head}%
%
%
\newcount\subheadno%
\global\subheadno=0%
\def\subheadinfo{\headinfo.\the\subheadno}%
\def\nextsubheadno{\global\advance\subheadno by 1 \global\procno=0 \global\rmkno=0 \subheadinfo}%
\makecounter{subhead}%
%
%
\newcount\procno%
\global\procno=0%
\def\procinfo{\subheadinfo.\the\procno}%
\def\nextprocno{\global\advance\procno by 1 \procinfo}%
\makecounter{proc}%
%
%
\newcount\figno%
\global\figno=0%
\def\figinfo{\subheadinfo.\the\figno}%
\def\nextfigno{\global\advance\figno by 1 \figinfo}%
\makecounter{fig}%
%
%
\newcount\eqnno%
\global\eqnno=0%
\def\eqninfo{\text{{\rm (\headinfo.\the\eqnno)}}}%
\def\nexteqnno[#1]{\global\advance\eqnno by 1 \eqninfo\hbox{\eqnlabel{#1}}}%
\makecounter{eqn}%
%
%
\newcount\rmkno%
\global\rmkno=0%
\def\rmkinfo{\text{\subheadinfo.\the\rmkno}}%
\def\nextrmkno[#1]{\global\advance\rmkno by 1 \rmkinfo\hbox{\rmklabel{#1}}}%
\makecounter{rmk}%
%
%
%
%
%
\def\gobbleeight#1#2#3#4#5#6#7#8{}%
\newcount\citationno%
\global\citationno=0%
\def\citationinfo{\the\citationno}%
\makecounter{citation}%
\newwrite\biblio%
\def\newref#1#2{%
\def\temptext{#2}%
\edef\bibliotextoutput{\expandafter\gobbleeight\meaning\temptext}%
\global\advance\citationno by 1\citationlabel{#1}%
\ifmakebiblio%
    \edef\fileoutput{\write\biblio{\noindent\hbox to 0pt{\hss$[\the\citationno]$}\hskip 0.2em\bibliotextoutput\medskip}}%
    \fileoutput%
\fi}%
\def\cite#1{%
$[\citationref{#1}]$%
\ifmakebiblio%
    \edef\fileoutput{\write\biblio{#1}}%
    \fileoutput%
\fi%
}%
%
%
%
%
\let\mypar=\par%
\edef\Pagetitle={Blank}\headline={\hfil\Pagetitle\hfil}%
\edef\Pagefooter={Blank}\footline={\hfil\Pagefooter\hfil}%
%
%
\newcount\showpagenumflag%
\global\showpagenumflag=0 %
\def\nextoddpage%
{\newpage\ifodd\pageno%
\else\global\showpagenumflag=0 %
\null\vfil\eject%
\global\showpagenumflag=1 %
\fi}%
%
%
\font\headfont=cmb12%
\def\newhead#1[#2]%
{\ifhmode\mypar\fi%
\ifnum\headno=0 \else\goodbreak\bigskip\fi%
{\headfont\noindent\nextheadno\ - #1.}\headlabel{#2}%
\nobreak\medskip}%
%
%
\def\newsubhead#1[#2]%
{\ifhmode\mypar\fi%
\ifnum\subheadno=0 \else\goodbreak\medskip\fi%
{\bf\noindent\nextsubheadno\ - #1.\ }\subheadlabel{#2}}%
%
%
\newif\ifinproclaim%
\global\inproclaimfalse%
\def\proclaim#1{%
\goodbreak\medskip
\bgroup\inproclaimtrue%
\noindent{\bf #1}%
\nobreak\medskip\sl}%
\def\noskipproclaim#1{%
\goodbreak\medskip%
\bgroup\inproclaimtrue%
\noindent{\bf #1}\nobreak\sl}%
\def\endproclaim{\mypar\egroup\nobreak\medskip\ignorespaces}%
%
%
%
\newcount\xpos\newcount\ypos
\def\makelabelgrid{%
\xpos=-5 \ypos=-5 %
\loop\ifnum\xpos<6 %
{\loop\ifnum\ypos<6 %
\def\labeltext{x}%
\ifnum\xpos=0\def\labeltext{+}\fi%
\ifnum\ypos=0\def\labeltext{+}\fi%
\placelabel[\xpos][\ypos]{\labeltext}%
\advance\ypos by 1 %
\repeat}%
\advance\xpos by 1 %
\repeat}%
\def\placelabel[#1][#2]#3{{%
\setbox10=\hbox{\raise #2cm \hbox{\hskip #1cm #3}}%
\ht10=0pt \dp10=0pt \wd10=0pt \box10}}%
%
%
%
%
\def\myitem#1{\noindent\hbox to .5cm{\hfill#1\hss}}%
%
%
%
%
%
%
%
%
%
\font\sansseriften=cmss10%
\font\sansserifseven=cmss7%
\font\sansseriffive=cmss5%
\newfam\sansseriffam%
\textfont\sansseriffam=\sansseriften%
\scriptfont\sansseriffam=\sansserifseven%
\scriptscriptfont\sansseriffam=\sansseriffive%
\def\mathsf{\fam\sansseriffam}%
%
%
%
\font\boldten=cmb10%
\font\boldseven=cmb7%
\font\boldfive=cmb5%
\newfam\mathboldfam%
\textfont\mathboldfam=\boldten%
\scriptfont\mathboldfam=\boldseven%
\scriptscriptfont\mathboldfam=\boldfive%
\def\mathbf{\fam\mathboldfam}%
%
%
%
\font\mycmmiten=cmmi10%
\font\mycmmiseven=cmmi7%
\font\mycmmifive=cmmi5%
\newfam\mycmmifam%
\textfont\mycmmifam=\mycmmiten%
\scriptfont\mycmmifam=\mycmmiseven%
\scriptscriptfont\mycmmifam=\mycmmifive%
\def\hexa#1{\ifcase #1 0\or 1\or 2\or 3\or 4\or 5\or 6\or 7\or 8\or 9\or A\or B\or C\or D\or E\or F\fi}%
\mathchardef\mathi="7\hexa\mycmmifam7B%
\mathchardef\mathj="7\hexa\mycmmifam7C%
%
%
\font\mymsbmten=msbm10 at 8pt%
\font\mymsbmseven=msbm7 at 5.6pt
\font\mymsbmfive=msbm5 at 4pt%
\newfam\mymsbmfam%
\textfont\mymsbmfam=\mymsbmten%
\scriptfont\mymsbmfam=\mymsbmseven%
\scriptscriptfont\mymsbmfam=\mymsbmfive%
\mathchardef\mybeth="7\hexa\mymsbmfam69%
\mathchardef\mygimmel="7\hexa\mymsbmfam6A%
\mathchardef\mydaleth="7\hexa\mymsbmfam6B%
%
%
%
%
\def\proof{{\noindent\bf Proof:\ }}%
\def\remark[#1]{{\noindent\bf Remark \nextrmkno[#1].}}%
\def\qed{~$\square$}%
\def\makeop#1{\global\expandafter\def\csname op#1\endcsname{{\text{#1}}}}%
\def\makeopsmall#1{\global\expandafter\def\csname op#1\endcsname{{\text{\lowercase{#1}}}}}%
%
%
\def\munion{\mathop{\cup}}%
\def\minter{\mathop{\cap}}%
%
%
\makeop{Ext}%
\makeop{Int}%
\makeop{Dist}%
\makeop{Diam}%
\makeop{Length}%
%
%
%
%
\def\minn{{m\in\Bbb{N}}}%
\def\mlimsup{\mathop{{\text{LimSup}}}}%
%
%
%
%
%
\makeop{Dim}%
\makeop{Ker}%
\makeop{Coker}%
\makeop{Tr}%
\makeop{Adj}%
\makeop{Det}%
\makeop{End}%
\makeop{Lin}%
\makeop{Symm}%
\makeop{Mult}%
%
%
\makeop{dx}%
\makeop{dy}%
\makeop{dz}%
\makeop{dt}%
\makeop{dVol}%
\makeop{dArea}%
\makeop{Supp}%
\makeop{Hess}%
\makeop{Lip}%
%
%
\makeop{Re}%
\makeop{Im}%
\makeop{Arg}%
\makeop{Log}%
\makeop{Exp}%
%
%
\makeopsmall{Cos}%
\makeopsmall{Sin}%
\makeopsmall{Tan}%
\makeopsmall{Sec}%
\makeopsmall{Cosec}%
\makeopsmall{Cot}%
\makeopsmall{ArcCos}%
\makeopsmall{ArcSin}%
\makeopsmall{ArcTan}%
\makeopsmall{ArcSec}%
\makeopsmall{ArcCosec}%
\makeopsmall{ArcCot}%
%
%
\makeopsmall{Cosh}%
\makeopsmall{Sinh}%
\makeopsmall{Tanh}%
\makeopsmall{ArcCosh}%
\makeopsmall{ArcSinh}%
\makeopsmall{ArcTanh}%
%
%
\makeop{Vol}%
\makeop{Area}%
\makeop{Riem}%
\makeop{Ric}%
\makeop{Scal}%
\makeop{Euc}%
\makeop{Imm}%
\makeop{Emb}%
%
%
\makeop{Id}%
\makeop{Ad}%
\makeop{O}%
\makeop{SO}%
\makeop{SL}%
\makeop{GL}%
\makeop{Conf}%
\makeop{Homeo}%
\makeop{Diff}%
\makeop{Isom}%
%
%
\makeop{Ind}%
\makeop{Sig}%
\makeop{Spec}%
%
%
\makeop{Conv}%
\makeop{Max}%
\makeop{Min}%
\makeop{Mod}%
\makeop{Deg}%
\makeop{loc}%
%
%
%
%
\def\Pagetitle{\hfil\ifnum\pageno=1\null\else{\rm $k$-surfaces in hyperbolic space.}\fi\hfil}
\def\Pagefooter{\hfil{\myfontdefault\folio}\hfil}
\catcode`\@=11
\def\triplealign#1{\null\,\vcenter{\openup1\jot \m@th %
\ialign{\strut\hfil$\displaystyle{##}\quad$&$\displaystyle{{}##}$\hfil&$\displaystyle{{}##}$\hfil\crcr#1\crcr}}\,}
\def\multiline#1{\null\,\vcenter{\openup1\jot \m@th %
\ialign{\strut$\displaystyle{##}$\hfil&$\displaystyle{{}##}$\hfil\crcr#1\crcr}}\,}
\catcode`\@=12
\makeop{I}%
\makeop{II}%
\makeop{III}%
\makeopsmall{Coth}%
\makeopsmall{Cossec}%
\makeop{dA}%
\makeop{PSL}%
\makeop{Hyp}%
\makeop{Ln}%
\makeop{A}%
\makeop{D}%
\makeop{J}%
\makeop{U}%
\makeop{T}%
\makeop{AdS}%
\makeop{sgn}%
\makeop{Spin}%
\makeop{Supp}%
\makeop{N}%
\def\opso{{\frak{so}}}
\makeop{T}%
\def\qi{{\mathbf{i}}}%
\def\qj{{\mathbf{j}}}%
\def\qk{{\mathbf{k}}}%
\def\Surface{S}%
\makeop{V}%
\makeop{H}%
\newref{Aron}{Aronszajn N., A unique continuation theorem for elliptic differential equations or inequalities of the second order, {\sl J. Math. Pures Appl.}, {\bf 36}, (1957), 235--239}
\newref{BarBegZegh}{Barbot T., B\'eguin F., Zeghib A., Prescribing Gauss curvature of surfaces in $3$-dimensional spacetimes: application to the Minkowski problem in the Minkowski space, {\sl Ann. Inst. Fourier}, {\bf 61}, no. 2, (2011), 511–-591}
\newref{BonMonSchI}{Bonsante F., Mondello G., Schlenker J. M., A cyclic extension of the earthquake flow I, {\sl Geom. Topol.}, {\bf 17}, no. 1, (2013), 157--234}
\newref{BonMonSchII}{Bonsante F., Mondello G., Schlenker J. M., A cyclic extension of the earthquake flow II, {\sl Ann. Sci. Ec. Norm. Sup\'er.}, {\bf 48}, no. 4, (2015), 811–-859}
\newref{Calabi}{Calabi E., Improper affine hyperspheres of convex type and a generalisation of a theorem by K. J\"orgens, {\sl Mich. Math. J.}, {\bf 5}, (1958), 105--126}
\newref{Girard}{Girard P. R., The quaternion group and modern physics, {\sl Eur. J. Phys.}, {\bf 5}, (1984) 25--32}
\newref{GromovPH}{Gromov P. H., Pseudo-holomorphic curves in symplectic manifolds, {\sl Inv. Math.}, {\bf 82}, (1985), 307--347}
\newref{Gromov}{Gromov M., Foliated Plateau problem, Part II: Harmonic maps of foliations, {\sl Geom. Func. Anal.}, {\bf 1}, no. 3, (1991), 253--320}
\newref{HarvLaws}{Harvey R., Lawson H. B. Jr., Calibrated geometries, {\sl Acta. Math.}, {\bf 148}, (1982), 47--157}
\newref{HarvLawsII}{Harvey F. R., Lawson H. B., Pseudoconvexity for the special Lagrangian potential equation, {\sl Calc. Var. PDEs.}, {\bf 60}, (2021)}
\newref{Jorgens}{J\"orgens K., \"Uber die L\"osungen der Differentialgleichung $rt-s^2=1$, {\sl Math. Ann.}, {\bf 127}, (1954), 130--134}
\newref{KokRossSajiUmeYam}{Kokubu M., Rossman W., Saji K., Umehara M., Yamada K., Singularities of flat fronts in hyperbolic space, {\sl Pac. J. Math.}, {\bf 221}, no. 2, (2005), 303--351}
\newref{KokRossUmeYam}{Kokubu M., Rossman W., Umehara M., Yamada K., Flat fronts in hyperbolic space and their caustics, {\sl J. Soc. Math. Japan}, {\bf 59}, no. 1, (2007), 265--299}
\newref{LabMP}{Labourie F., M\'etriques prescrites sur le bord des vari\'et\'es hyperboliques de dimension $3$, {\sl J. Diff. Geom.}, {\bf 35}, no. 3, (1992), 609–-626}
\newref{LabMA}{Labourie F., Probl\`emes de Monge-Amp\`ere, courbes pseudo-holomorphes et laminations, {\sl Geom. Func. Anal.}, {\bf 7}, (1997), 496--534}
\newref{LychRub}{Lychagin V. V., Rubtsov V. N., Local classification of Monge-Amp\`ere differential equations (Russian), {\sl Dokl. Akad. Nauk SSSR}, {\bf 272}, no. 1, (1983), 34--38}
\newref{LychRubChek}{Lychagin V. V., Rubtsov V. N., Chekalov I. V., A classification of Monge-Amp\`ere equations, {\sl Ann. Sci. \'Ecole Norm. Sup.}, {\bf 26}, no. 3, (1993), 281--308}
\newref{MartMil}{Mart\'\i nez A., Mil\'an F., Flat fronts in hyperbolic $3$-space with prescribed singularities, {\sl Ann. Glob. Anal. Geom.}, {\bf 46}, (2014), 227--239}
\newref{McDuffSal}{McDuff D., Salamon D., {\sl J-holomorphic curves and quantum cohomology}, University Lecture Series, {\bf 6}, AMS, Providence, (1994)}
\newref{Pogorelov}{Pogorelov A. V., On the improper convex affine hyperspheres, {\sl Geom. Dedi.}, {\bf 1}, (1972), 33--46}
\newref{Rub}{Rubtsov V., Geometry of Monge-Ampère structures, in {\sl Nonlinear PDEs, their geometry, and applications}, Birkh\"auser/Springer, (2019), 95--156}
\newref{Sch}{Schlenker J. M., Hyperbolic manifolds with convex boundary, {\sl Inv. Math.}, {\bf 163}, (2006), 109--169}
\newref{Spivak}{Spivak M., {\sl A comprehensive introduction to differential geometry}, Vol. IV, Publish or Perish, (1999)}
\newref{SmiAA}{Smith G., An Arzela-Ascoli Theorem for Immersed Submanifolds, {\sl Ann. Fac. Sci. Toulouse Math.}, {\bf 16}, no. 4, (2007), 817--866}
\newref{SmiSLC}{Smith G., Special Lagrangian curvature, {\sl Math. Ann.}, {\bf 335}, no. 1, (2013), 57--95}
\newref{SmiAS}{Smith G., On the asymptotic geometry of finite-type k-surfaces in three-dimensional hyperbolic space, arXiv:1908.04834}
\newref{TouLabWol}{Toulisse J., Labourie F., Plateau Problems for Maximal Surfaces in Pseudo-Hyperbolic Spaces, to appear in {\sl Ann. Sci. \'Ec. Norm. Sup\'er}, arXiv:2006.12190}
\def\centre{\rightskip=0pt plus 1fil \leftskip=0pt plus 1fil \spaceskip=.3333em \xspaceskip=.5em \parfillskip=0em \parindent=0em}%
\def\textmonth#1{\ifcase#1\or January\or Febuary\or March\or April\or May\or June\or July\or August\or September\or October\or November\or December\fi}
\font\abstracttitlefont=cmr10 at 14pt {\abstracttitlefont\centre Quaternions, Monge--Amp\`ere structures and $k$-surfaces\par}
\bigskip
{\centre 5th October 2022\par}
\bigskip
{\centre Graham Smith\par}
\bigskip
\noindent{\bf Abstract:~}In \cite{LabMA} Labourie develops a theory of immersed surfaces of prescribed extrinsic curvature which has since found widespread applications in hyperbolic geometry, general relativity, Teichm\"uller theory, and so on. In this chapter, we present a quaternionic reformulation of these ideas. This yields simpler proofs of the main results whilst pointing towards the higher-dimensional generalisation studied by the author in \cite{SmiSLC}.
\bigskip
\noindent{\bf Classification AMS~:~}53A05, 12E15, 35J96
\bigskip
%
%
\newhead{Introduction}[Introduction]
\newsubhead{Introduction}[Introduction]
Immersed surfaces of prescribed extrinsic curvature in $3$-dimensional manifolds have fascinated mathematicians for almost two centuries. Following the pioneering work \cite{LabMA} of Labourie, remarkable developments have been made in our understanding of these objects, leading to striking applications across a broad range of mathematical theories\numberedfootnote{The reader may consult, for example, \cite{BarBegZegh}, \cite{BonMonSchI}, \cite{BonMonSchII}, \cite{LabMP}, \cite{Sch}, \cite{SmiAS}, and so on, for a selection of applications of these techniques.}. In this chapter, we propose a quaternionic reformulation of Labourie's ideas. Not only will this yield simpler proofs of the main results, but it will also point towards their higher-dimensional generalisations studied by the author in \cite{SmiSLC}.
\par
We first recall the main elements of Labourie's work as it pertains to immersed surfaces. Let $X:=(X,h)$ be a complete, oriented, $3$-dimensional riemannian manifold, let $TX$ denote its tangent bundle, and let $SX\subseteq TX$ denote its unit sphere bundle. We define an (oriented) {\bf immersed surface} in $X$ to be a pair $(\Surface,e)$, where $\Surface$ is an oriented surface, and $e:\Surface\rightarrow X$ is a smooth immersion. Given such a pair, we denote by $\nu_e:\Surface\rightarrow SX$ its unit normal vector field compatible with the orientation, by $\opI_e$, $\opII_e$ and $\opIII_e$ its first, second and third fundamental forms respectively, by $A_e$ its shape operator, and by $K_e:=\opDet(A_e)$ its {\bf extrinsic curvature function}. We say that the immersed surface is {\bf infinitesimally strictly convex (ISC)} whenever its second fundamental form is positive definite, we say that it is {\bf quasicomplete} whenever it is complete with respect to the riemannian metric $\opI_e+\opIII_e$, and, given a smooth function $\kappa:SX\rightarrow\Bbb{R}$, we say that its extrinsic curvature is {\bf prescribed} by $\kappa$ whenever
$$
K_e := \kappa\circ\nu_e.\eqnum{\nexteqnno[CurvaturePrescription]}
$$
\par
We denote $\hat{e}:=\nu_e$, and we call the immersed surface $(\Surface,\hat{e})$ the {\bf Gauss lift} of $(\Surface,e)$. Note that quasicompleteness of $(\Surface,e)$ is equivalent to completeness of its Gauss lift. Labourie's key insight is that the Gauss lift of any ISC immersed surface of prescribed extrinsic curvature is a pseudo-holomorphic curve for some suitable almost complex structure. This allows the powerful theory developed by Gromov in \cite{GromovPH} to be applied. The first consequence is the following compactness result.
\proclaim{Theorem \nextprocno, {\bf Labourie's compactness theorem}}
\noindent Let $\kappa:SX\rightarrow\Bbb{R}$ be a smooth, positive function, and let $(\Surface_m,e_m,p_m)$ be a sequence of quasicomplete, pointed, ISC immersed surfaces in $X$ of extrinsic curvature prescribed by $\kappa$. If the sequence $(e_m(p_m))$ is precompact in $X$, then the sequence $(\Surface_m,\hat{e}_m,p_m)$ of Gauss lifts is precompact in the smooth Cheeger--Gromov topology.
\endproclaim
\proclabel{IntroLabouriesCompactnessTheorem}
\remark[IntroLabouriesCompactnessTheoremI] This is proven in Theorem \procref{LabouriesCompactnessTheoremPGCI}. The smooth Cheeger--Gromov topology is described in Section \subheadref{Compactness}.
\medskip
\remark[IntroLabouriesCompactnessTheoremII] We have stated this result in its simplest possible form. It may however be generalised in a number of ways. For example, the ambient space $X$ as well as the function $\kappa$ can be allowed to vary, the conditions of completeness of the ambient space and quasicompleteness of the immersed surfaces can also be relaxed, and so on. We refer the reader to Section \subheadref{Compactness}, where potential generalisations are explained in greater detail. The reader may likewise consult \cite{TouLabWol} for recent developments concerning the case with non-trivial boundary.
\medskip
\noindent Significantly, Theorem \procref{IntroLabouriesCompactnessTheorem} is only of limited use without an understanding of degenerate limits, that is, those limits that are not Gauss lifts of immersed surfaces. This is addressed by Labourie in his second key result. For any complete geodesic $\Gamma$ in $X$, we denote by $\opN\Gamma$ its unit normal bundle. We define a {\bf tube} in $SX$ to be an immersed surface $(\Surface,\hat{e})$ which is a cover of $\opN\Gamma$ for some complete geodesic $\Gamma$.
\proclaim{Theorem \nextprocno, {\bf Labourie's dichotomy}}
\noindent With the notation of Theorem \procref{IntroLabouriesCompactnessTheorem}, every accumulation point of the sequence $(\Surface_m,e_m,p_m)$ is either a tube or the Gauss lift of some quasicomplete, pointed, ISC immersed surface in $X$ of extrinsic curvature prescribed by $\kappa$.
\endproclaim
\proclabel{IntroLabouriesDichotomy}
\remark[IntroLabouriesDichotomyI] This is proven in Theorem \procref{LabouriesCompactnessTheoremPGCI}.
\medskip
\noindent These two remarkable results and their variants have since become the basis of a rich theory of surfaces of prescribed extrinsic curvature in riemannian and semi-riemannian manifolds. It will be the object of this chapter to present a quaternionic framework within which they can be proved, allowing us, on the one hand, to emphasize their hyperk\"ahler nature, and, on the other, to point towards their higher-dimensional generalisations.
\newsubhead{Acknowledgements}[Acknowledgements]
\noindent This paper was in part written whilst the author was visiting the Institut des Hautes \'Etudes Scientifiques. The author is grateful for the excellent working conditions enjoyed during that stay.
\newhead{Quaternions and Bernstein-type theorems}[QuaternionsAndBernsteinTypeTheorems]
\newsubhead{Quaternions}[Quaternions]
We begin our work with a detailed review of Hamilton's theory of quaternions. Introduced in 1843, this theory has enjoyed a striking revival over the past forty years on account of the remarkably simple approaches it provides to various deep mathematical and physical phenomena (see, for example, \cite{Girard}). Its application to the study of prescribed curvature surfaces, which will be addressed in detail in the sequel, presents yet another instance of the ubiquity of this theory that makes it so intriguing.
\par
Let $\Bbb{H}$ denote the algebra of quaternions. This is the associative, unital algebra over $\Bbb{R}$ generated by the $3$ elements $\qi$, $\qj$ and $\qk$, with the relations
$$
\qi^2 = \qj^2 = \qk^2 = \qi\cdot\qj\cdot\qk = -1.\eqnum{\nexteqnno[QuaternionRelations]}
$$
Given an element $x\in\Bbb{H}$ of the form
$$
x:=a + b\qi + c\qj + d\qk,\eqnum{\nexteqnno[GeneralQuaternionFormula]}
$$
its {\bf conjugate} is defined by
$$
\overline{x}:=a-b\qi-c\qj-d\qk.\eqnum{\nexteqnno[DefinitionOfQuaternionConjugate]}
$$
Conjugation is an anti-involution of $\Bbb{H}$ in the sense that, for all $x,y$,
$$
\overline{x\cdot y} = \overline{y}\cdot\overline{x}.\eqnum{\nexteqnno[ComplexConjugationsAnticommutes]}
$$
Its $(+1)$-eigenspace, the space of {\bf real quaternions}, is a $1$-dimensional subspace which we denote by $\Cal{R}$. Its $(-1)$-eigenspace, the space of {\bf imaginary quaternions}, is a $3$-dimensional subspace which we denote by $\Cal{I}$. We likewise denote by $\Cal{R}$ and $\Cal{I}$ the respective projections onto these subspaces.
\par
By \eqnref{ComplexConjugationsAnticommutes}, for all $x\in\Bbb{H}$,
$$
\overline{x\cdot\overline{x}} = \overline{\overline{x}}\cdot\overline{x} = x\cdot\overline{x},\eqnum{\nexteqnno[RealityOfInnerProduct]}
$$
so that $x\cdot\overline{x}$ is always real. We thus define an inner product by
$$
\langle x,y\rangle := \Cal{R}(x\cdot\overline{y}).\eqnum{\nexteqnno[QuaternionInnerProduct]}
$$
With $x$ as in \eqnref{GeneralQuaternionFormula},
$$
\|x\|^2 = a^2 + b^2 + c^2 + d^2,\eqnum{\nexteqnno[FormulaForNormSquared]}
$$
so that \eqnref{QuaternionInnerProduct} coincides with the standard inner product of $4$-dimensional euclidian space. Furthermore, for all $x,y\in\Bbb{H}$,
$$
\|x\cdot y\|^2 = x\cdot y\cdot\overline{x\cdot y} = x\cdot y\cdot\overline{y}\cdot\overline{x} = \|x\|^2\|y\|^2,\eqnum{\nexteqnno[LengthIsMultiplicative]}
$$
so that length is multiplicative.
\par
For all $x,y,z\in\Cal{I}$,
$$
\Cal{R}(x\cdot y\cdot z) = \Cal{R}(\overline{x}\cdot\overline{y}\cdot z) = \Cal{R}(\overline{y\cdot x}\cdot z) = \Cal{R}(y\cdot x\cdot\overline{z}) = -\Cal{R}(y\cdot x\cdot z).\eqnum{\nexteqnno[AntisymmetryI]}
$$
We thus define an alternating $3$-form over $\Cal{I}$ by
$$
\omega(x,y,z) := -\Cal{R}(x\cdot y\cdot z).\eqnum{\nexteqnno[VolumeForm]}
$$
This is the volume form of $\langle\cdot,\cdot\rangle$ with orientation chosen such that $(\qi,\qj,\qk)$ is a positive triple.
\par
By \eqnref{LengthIsMultiplicative}, the sphere $\Bbb{S}^3$ of unit quaternions is a subgroup of $\Bbb{H}$ with inverse given by conjugation.
\proclaim{Lemma \& Definition \nextprocno}
\noindent The homomorphism $h:\Bbb{S}^3\times\Bbb{S}^3\rightarrow\opSO(\Bbb{H})$ given by
$$
h(x,y)z := x\cdot z\cdot \overline{y}.\eqnum{\nexteqnno[CoveringOfSO]}
$$
is a double cover of $\opSO(\Bbb{H})$. In particular, it identifies $\Bbb{S}^3\times\Bbb{S}^3$ with $\opSpin(\Bbb{H})$.
\endproclaim
\proclabel{CoveringOfSO}
\proof Since it is a continuous homomorphism between Lie groups of the same dimension, it suffices to show that its kernel is
$$
\opKer(h) = \left\{ (1,1),(-1,-1)\right\}.\eqnum{\nexteqnno[KernelOfH]}
$$
However, suppose that $h(x,y)=\opId$. Substituting $z=1$ into \eqnref{CoveringOfSO} yields $x\cdot\overline{y}=1$, so that $x=y$. Thus, for all $z$,
$$
z\cdot x = (h(x,x)\cdot z)\cdot x = x\cdot z\cdot\overline{x}\cdot x = x\cdot z,
$$
so that
$$
[x,z] = 0.
$$
Since $z$ is arbitrary, $x$ is real, and since $x$ has unit norm, it is equal $\pm 1$. This proves \eqnref{KernelOfH}, and the result follows.\qed
\newsubhead{Compatible complex structures}[CompatibleComplexStructures]
Recall that a {\bf complex structure} over $\Bbb{H}$ is an $\Bbb{R}$-linear map $J:\Bbb{H}\rightarrow\Bbb{H}$ such that
$$
J^2 = -\opId.\eqnum{\nexteqnno[ComplexStructureCondition]}
$$
We say, in addition, that $J$ is {\bf compatible} whenever it preserves the metric, that is, whenever, for all $x\in\Bbb{H}$,
$$
\|Jx\| = \|x\|.\eqnum{\nexteqnno[JPreservesLength]}
$$
We now proceed to identify all compatible complex structures over $\Bbb{H}$.
\proclaim{Lemma \nextprocno}
\noindent The set of square roots of $-1$ in $\Bbb{H}$ is the sphere $\Bbb{S}^3\minter\Cal{I}$ of unit, imaginary quaternions.
\endproclaim
\proclabel{ComplexStructures}
\proof Indeed, $x^2=-1$ if and only if $x^{-1}=-x$. This holds if and only if $\|x\|=1$ and $x=-\overline{x}$, as desired.\qed
\medskip
\noindent For every unit, imaginary quaternion $x$, $h(x,1)$ and $h(1,x)$ are trivially compatible complex structures over $\Bbb{H}$. We now verify that there are no others.
\proclaim{Lemma \& Definition \nextprocno}
\noindent If $J:\Bbb{H}\rightarrow\Bbb{H}$ is a compatible complex structure, then $J$ is given by multiplication either on the left or on the right by a unit, imaginary quaternion. We call the former {\bf left complex structures} and the latter {\bf right complex structures}.
\endproclaim
\proclabel{LeftAndRightComplexStructures}
\remark[LeftAndRightComplexStructures] Every left complex structure trivially commutes with every right complex structure. Note also that quaternionic conjugation sends left complex structures into right complex structures and vice-versa.
\medskip
\proof Since $J\in\opSO(\Bbb{H})$, by Lemma \procref{CoveringOfSO}, there exists $(x,y)\in\Bbb{S}^3\times\Bbb{S}^3$ such that $J=h(x,y)$. Since
$$
-J^2=h(x^2,-y^2)=\opId,
$$
it follows that
$$
x^2=-y^2=\pm 1.
$$
Without loss of generality, we may suppose that $(x^2,y^2)=(1,-1)$. It then follows that $x=\pm1$, $y\in\Bbb{S}^3\minter\Cal{I}$, and
$$
J = h(\pm1,y) = h(1,\pm y),
$$
as desired.\qed
\medskip
\noindent In particular, we obtain the following algebraic characterisation of left and right complex structures.
\proclaim{Lemma \nextprocno}
\noindent Let $J:\Bbb{H}\rightarrow\Bbb{H}$ be a right complex structure. If $J':\Bbb{H}\rightarrow\Bbb{H}$ is another compatible complex structure which commutes with $J$, then either $J'$ is a left complex structure, or $J'=\pm J$.
\endproclaim
\proclabel{CharacterisationOfLeftComplexStructures}
\remark[CharacterisationOfLeftComplexStructures] An analogous result trivially holds with the roles of left and right complex structures inverted.
\medskip
\proof It suffices to show that if $J'$ is a right complex structure which commutes with $J$, then $J'=\pm J$. However, let $x$ and $y$ be unit, imaginary quaternions such that $J=h(1,x)$ and $J'=h(1,y)$. Since $J$ and $J'$ commute, so too do $x$ and $y$. Thus, since $x$ and $y$ are imaginary,
$$
\overline{x\cdot\overline{y}} = y\cdot\overline{x} = -y\cdot x = -x\cdot y = x\cdot\overline{y},
$$
so that $x\cdot\overline{y}$ is real. Since $x$ and $y$ both have unit length, this holds if and only if $x=\pm y$, and the result follows.\qed
\newsubhead{Compatible quaternionic structures}[HomeomorphismsAndCompatibleQuaternionicStructures]
Let $E:=(E,g)$ be an inner-product space. We define a {\bf compatible quaternionic structure} over $E$ to be an algebra homomorphism $\rho:\Bbb{H}\rightarrow\opEnd(E)$ such that, for all $x\in\Bbb{H}$ and for all $u\in E$,
$$
\|\rho(x)u\| = \|x\|\cdot\|u\|.\eqnum{\nexteqnno[CompatibilityCondition]}
$$
We now proceed to identify all compatible quaternionic structures over $\Bbb{H}$.
\proclaim{Lemma \nextprocno}
\noindent If $\alpha:\Bbb{H}\rightarrow\Bbb{H}$ be an non-trivial algebra homomorphism, then there exists a unit quaternion $z$ such that, for all $x\in\Bbb{H}$,
$$
\alpha(x) = z\cdot x\cdot \overline{z}.\eqnum{\nexteqnno[HomomorphismsOfQuaternions]}
$$
\endproclaim
\proclabel{HomomorphismsOfQuaternions}
\proof Since $\Bbb{H}$ is a skew field, $\alpha$ is injective and $\alpha(\pm 1)=\pm 1$. In particular, $\alpha$ preserves $\Cal{R}$. It likewise preserves the set $\Bbb{S}^3\minter\Cal{I}$ of square roots of $-1$, and therefore also $\Cal{I}$. It follows that $\alpha$ preserves the decomposition $\Bbb{H}=\Cal{R}\oplus\Cal{I}$. It therefore also preserves conjugation, and thus also the metric \eqnref{QuaternionInnerProduct}. Finally, by \eqnref{VolumeForm}, $\alpha$ preserves the orientation of $\Cal{I}$ and thus also of $\Bbb{H}$, and is consequently an element of $\opSO(\Bbb{H})$. It follows by Lemma \procref{CoveringOfSO} that there exist unit quaternions $z,w\in\Bbb{S}^3$ such that $\alpha=h(w,z)$. Finally, since $\alpha(1)=1$, $w=z$, and the result follows.\qed
\proclaim{Lemma \& Definition \nextprocno}
\noindent If $\alpha:\Bbb{H}\rightarrow\opEnd(\Bbb{H})$ is a compatible quaternionic structure, then either
\medskip
\myitem{(1)} there exists $z\in\Bbb{S}^3$ such that, for all $x,y\in\Bbb{H}$,
$$
\alpha(x)\cdot y = z\cdot x\cdot\overline{z}\cdot y,\ \text{or}\eqnum{\nexteqnno[LeftQuaternionicStructure]}
$$
\myitem{(2)} there exists $z\in\Bbb{S}^3$ such that, for all $x,y\in\Bbb{H}$,
$$
\alpha(x)\cdot y = y\cdot z\cdot\overline{x}\cdot\overline{z}.\eqnum{\nexteqnno[RightQuaternionicStructure]}
$$
\noindent We call the former {\bf left quaternionic structures} and the latter {\bf right quaternionic structures}.
\endproclaim
\proclabel{LeftAndRightQuaternionicStructures}
\remark[LeftAndRightQuaternionicStructures] As before, every left quaternionic structure commutes with every right quaternionic structure. Furthermore, quaternionic conjugation sends left quaternionic structures to right quaternionic structures and vice-versa.
\medskip
\proof Since $\Bbb{H}$ is a skew field, $\alpha$ is injective and $\alpha(\pm 1)=\pm\opId$. It follows that $\alpha$ maps the set of unit, imaginary quaternions to the set of compatible complex structures of $\Bbb{H}$. By connectedness, we may suppose without loss of generality that, for every unit, imaginary quaternion $x$, $\alpha(x)$ is a left complex structure. It follows that, for all $x,y\in\Bbb{H}$,
$$
\alpha(x)\cdot y = \rho(x)\cdot y,
$$
for some non-trivial algebra homomorphism $\rho:\Bbb{H}\rightarrow\Bbb{H}$. The result now follows by Lemma \procref{HomomorphismsOfQuaternions}.\qed
\newsubhead{Calibrations}[Calibrations]
The utility of quaternions to the theory of partial differential equations arises from the theory of calibrations developed by Harvey \& Lawson in \cite{HarvLaws}. We now explain how this applies our setting. For every unit, imaginary quaternion $x$, we define
$$
J_x := h(x,1),\eqnum{\nexteqnno[LeftComplexStructure]}
$$
and we define the symplectic form $\omega_x$ over $\Bbb{H}$ by
$$
\omega_x := \langle\cdot,J_x\cdot\rangle.\eqnum{\nexteqnno[AssociatedSymplecticForm]}
$$
\proclaim{Lemma \nextprocno}
\noindent Let $(x,y,z)$ be an orthonormal triplet of unit, imaginary quaternions. Let $P\subseteq\Bbb{H}$ be a real plane. The area form $\opdA$ of $P$ satisfies, for all $\xi,\nu\in P$,
$$
dA(\xi,\nu)^2 = \omega_x(\xi,\nu)^2 + \omega_y(\xi,\nu)^2 + \omega_z(\xi,\nu)^2.\eqnum{\nexteqnno[Calibration]}
$$
\endproclaim
\proclabel{CalibrationA}
\remark[CalibrationA] In other words, the triplet $(\omega_x,\omega_y,\omega_z)$ forms a {\bf calibration} of $\Bbb{H}$ in the sense of \cite{HarvLaws}.
\medskip
\proof We may suppose that $\xi$ and $\nu$ both have unit length. Since $(\nu,J_x\nu,J_y\nu,J_z\nu)$ is an orthonormal real basis of $\Bbb{H}$,
$$
\langle\xi,\nu\rangle^2 + \langle\xi,\opJ_x\nu\rangle^2 + \langle\xi,\opJ_y\nu\rangle^2 + \langle\xi,\opJ_z\nu\rangle^2 = \|\xi\|^2 = 1.
$$
Thus, if $\theta$ denotes the angle between $\xi$ and $\nu$, then
$$
dA(\xi,\nu)^2 = \opSin^2(\theta) = 1-\langle\xi,\nu\rangle^2 = \omega_x(\xi,\nu)^2 + \omega_y(\xi,\nu)^2 + \omega_z(\xi,\nu)^2,
$$
as desired.\qed
\medskip
\noindent In particular, we obtain the following result for real planes in $\Bbb{H}$.
\proclaim{Lemma \nextprocno}
\noindent Let $x$ be a unit, imaginary quaternion. Let $P$ be a real plane in $\Bbb{H}$. $P$ is $\opJ_x$-complex if and only if it is $\omega_y$-lagrangian for every unit, imaginary quaternion $y$ orthogonal to $x$.
\endproclaim
\proclabel{CalibrationB}
\proof Let $y$ and $z$ be unit, imaginary quaternions such that $(x,y,z)$ is an orthonormal triplet, and let $\opdA$ denote the area form of $P$. $P$ is $J_x$-complex if and only if, for all $\xi,\nu\in P$,
$$
\opdA(\xi,\nu)^2 = \langle\xi,\opJ_x\nu\rangle^2 = \omega_x(\xi,\nu)^2.
$$
By \eqnref{Calibration}, this holds if and only if $\omega_y$ and $\omega_z$ both vanish over $P$, as desired.\qed
\newsubhead{The Monge--Amp\`ere equation}[TheMongeAmpereEquation]
We now apply our abstract framework to the study of solutions of the real Monge--Amp\`ere equation
$$
\opDet(\opHess(u)) = 1.\eqnum{\nexteqnno[MongeAmpere]}
$$
To this end, we identify $\Bbb{C}$ with $\Bbb{R}^2$ in the natural manner, and we introduce an explicit compatible quaternionic structure over $\Bbb{C}\oplus\Bbb{C}=\Bbb{R}^2\oplus\Bbb{R}^2$ as follows. First, let $J_0$ denote the operator of multiplication by $i$, so that
$$
J_0 := \pmatrix 0\hfill& -1\hfill\cr 1\hfill& 0\hfill\cr\endpmatrix,\eqnum{\nexteqnno[StandardComplexStructure]}
$$
and let $g$ denote the standard metric over $\Bbb{C}\oplus\Bbb{C}$, that is
$$
g((z,w)^t,(z,w)^t) = \left|z\right|^2 + \left|w\right|^2.\eqnum{\nexteqnno[StandardMetric]}
$$
We define
$$\eqalign{
I &:= \pmatrix J_0\hfill&0\hfill\cr 0\hfill&-J_0\hfill\cr\endpmatrix,\cr
J &:= \pmatrix 0\hfill&J_0\hfill\cr J_0\hfill&0\hfill\cr\endpmatrix,\ \text{and}\cr
K &:= \pmatrix 0\hfill&-\opId\hfill\cr \opId\hfill&0\hfill\cr\endpmatrix.\cr}\eqnum{\nexteqnno[QuaternionicStructure]}
$$
Note that, in the spirit of Lemma \procref{CharacterisationOfLeftComplexStructures}, $I$, $J$ and $K$ are generators of the set of compatible complex structures over $\Bbb{C}\oplus\Bbb{C}$ which commute with
$$
\hat{J}_0 := \pmatrix J_0\hfill&0\hfill\cr 0\hfill&J_0\hfill\cr\endpmatrix.\eqnum{\nexteqnno[BaseComplexStructure]}
$$
Let $(\omega_i,\omega_j,\omega_k)$ denote the triplet of symplectic forms corresponding to $(I,J,K)$, that is
$$\eqalign{
\omega_i&:=g(\cdot,I\cdot),\cr
\omega_j&:=g(\cdot,J\cdot),\ \text{and}\cr
\omega_k&:=g(\cdot,K\cdot).\cr}\eqnum{\nexteqnno[SymplecticStructures]}
$$
\par
We now examine the algebraic properties of different types of lagrangian subspaces of $\Bbb{C}\oplus\Bbb{C}$. To this end, let $P\subseteq\Bbb{C}\oplus\Bbb{C}$ be a real plane which is a graph over the first component, that is
$$
P:=\left\{(x,Ax)\ |\ x\in\Bbb{R}^2\right\},\eqnum{\nexteqnno[GraphPlane]}
$$
for some matrix $A\in\opEnd(2)$. The following two useful relations for matrices $A\in\opEnd(2)$ may be verified by inspection.
$$\eqalign{
A^t J_0 A J_0 &= -\opDet(A)\opId,\ \text{and}\cr
A-A^t &= -\opTr(A J_0)J_0.\cr}\eqnum{\nexteqnno[TwoDimMatrixRelations]}
$$
\proclaim{Lemma \nextprocno}
\myitem{(1)} $P$ is $\omega_i$-lagrangian if and only if $\opDet(A)=1$;
\medskip
\myitem{(2)} $P$ is $\omega_j$-lagrangian if and only if $\opTr(A)=0$; and
\medskip
\myitem{(3)} $P$ is $\omega_k$-lagrangian if and only if $A=A^t$.
\endproclaim
\proclabel{FirstLexicon}
\proof Indeed, $P$ is $\omega_i$-lagrangian if and only if,
$$\triplealign{
&\langle x,J_0y\rangle - \langle Ax,J_0Ay\rangle &=0\ \forall x,y\in\Bbb{R}^2\cr
\Leftrightarrow&\langle x,(\opId+A^t J_0 A J_0)y\rangle &=0\ \forall x,y\in\Bbb{R}^2.\cr}
$$
This holds if and only if
$$\triplealign{
\Leftrightarrow&\opId + A^t J_0 A J_0 &=0\cr
\Leftrightarrow&\opDet(A) &= 1,\cr}
$$
where the last equivalence follows by \eqnref{TwoDimMatrixRelations}. In a similar manner, using \eqnref{TwoDimMatrixRelations} again, we show that $P$ is $\omega_j$-lagrangian if and only if $\opTr(A)=0$. Finally, it is a standard result that $P$ is $\omega_k$-lagrangian if and only if $A=A^t$, and this completes the proof.\qed
\medskip
\noindent We thus obtain the following key result which relates solutions of the $2$-dimensional Monge--Amp\`ere equation to pseudo-holomorphic curves.
\proclaim{Lemma \nextprocno, {\bf Monge--Amp\`ere\ =\ Pseudoholomorphic}}
\noindent Let $\Omega$ be a simply-connected open subset of $\Bbb{R}^2$. Let $\alpha:\Omega\rightarrow\Bbb{R}^2$ be a smooth function and let $\Surface\subseteq\Omega\times\Bbb{R}^2$ denote its graph. Then $\alpha$ is the derivative of a smooth solution $u:\Omega\rightarrow\Bbb{R}$ of the real Monge--Amp\`ere equation \eqnref{MongeAmpere} if and only if $S$ is a $J$-holomorphic curve.
\endproclaim
\proclabel{SolutionsOfMAAreJHolomorphicCurves}
\proof Indeed, by Lemma \procref{CalibrationB}, $\Surface$ is $J$-holomorphic if and only if it is $\omega_i$- and $\omega_k$-lagrangian. By Lemma \procref{FirstLexicon}, this holds if and only if $\opDet(D\alpha)=1$ and $\alpha$ is closed. Since $\Omega$ is simply-connected, $\alpha$ is closed if and only if it is the derivative of some smooth function $u:\Omega\rightarrow\Bbb{R}$. The result now follows, since $\opDet(\opHess(u))=\opDet(D\alpha)=\pm 1$.\qed
\newsubhead{Positivity I}[Positivity]
We now express convexity in the quaternionic framework. To this end, let $\pi_1,\pi_2:\Bbb{C}\oplus\Bbb{C}\rightarrow\Bbb{C}$ denote respectively the projections onto the first and second components, and denote
$$
V:=\opKer(\pi_1)=\left\{0\right\}\oplus\Bbb{C}.\eqnum{\nexteqnno[VerticalSubspace]}
$$
We call $V$ the {\bf vertical subspace} of $\Bbb{C}\oplus\Bbb{C}$. Let $m$ be the symmetric bilinear form defined over $\Bbb{C}\oplus\Bbb{C}$ by
$$
m((z_1,w_1)^t,(z_2,w_2)^t) := \langle z_1,w_2\rangle + \langle z_2,w_1\rangle.\eqnum{\nexteqnno[MinkowskiMetric]}
$$
Note that $m$ vanishes over $V$ and that
$$
m(\cdot,\cdot) = -m(I\cdot,I\cdot) = m(J\cdot,J\cdot) = -m(K\cdot,K\cdot).\eqnum{\nexteqnno[InvarianceOfM]}
$$
It is worth noting that $m$ is, up to sign, uniquely determined by these properties.
\proclaim{Lemma \nextprocno}
\noindent Up to sign, $m$ is the unique symmetric bilinear form of unit norm satisfying \eqnref{InvarianceOfM} which vanishes over $V$.
\endproclaim
\proclabel{CharacterisationOfM}
\proof Indeed, let $m'$ be another symmetric bilinear form which vanishes over $V$ and which satisfies \eqnref{InvarianceOfM}. Since $m'$ vanishes over $V$ and is $J$-invariant, it is uniquely determined by the restriction of $\tilde{m}':=m(\cdot,J\cdot)$ to $V$. However, since $m'$ is $I$-anti-invariant, and since $I$ anticommutes with $J$, $\tilde{m}'$ is also $I$-invariant, and its restriction to $V$ is thus unique up to a scalar factor. Finally, since $m'$ has unit norm, it follows that $m'=\pm m$, as desired.\qed
\medskip
\noindent Define the subgroup $\Sigma\subseteq\opSO(4)$ by
$$
\Sigma := \left\{ \pmatrix A\hfill& 0\hfill\cr 0\hfill& A\hfill\cr\endpmatrix\ \bigg|\ A\in\opSO(2)\right\}.\eqnum{\nexteqnno[StabiliserSubgroup]}
$$
\proclaim{Lemma \nextprocno}
\noindent $\Sigma$ is the stabiliser subgroup in $\opSO(4)$ both of $(I,J,K,m)$ and of $(J,V,\omega)$, where $\omega$ as any orientation form of $V$.
\endproclaim
\proclabel{StabiliserSubgroup}
\remark[StabiliserSubgroup] In particular, the prescription of a triplet of the form $(V,\omega,J)$ over any $4$-dimensional inner-product space $(E,g)$ is equivalent to the prescription of a compatible quaternionic structure $(I,J,K)$ together with a symmetric bilinear form $m$ satisfying \eqnref{InvarianceOfM}.
\medskip
\proof Suppose that $M\in\opSO(4)$ preserves $(I,J,K,m)$. Define the involution $\alpha$ of $\Bbb{C}\oplus\Bbb{C}$ by $\alpha(z,w)= (z,\overline{w})$, and denote $\tilde{M}:=\alpha M\alpha$. Since $(\alpha I\alpha)(z,w)=(iz,iw)$, $\tilde{M}\in\opU(2)$. Since, in addition, $\tilde{M}$ preserves $(\alpha K\alpha)(z,w)=(-\overline{w},\overline{z})$,
$$
M = \pmatrix a\hfill& b\hfill\cr -\overline{b}\hfill&\overline{a}\hfill\cr\endpmatrix,
$$
for some $a,b\in\Bbb{C}$ such that $\left|a\right|^2+\left|b\right|^2=1$. Finally, since $\tilde{M}$ preserves $m$, $b=0$, and it follows that $M$ is an element of $\Sigma$, as desired.
\par
Suppose now that $M\in\opSO(4)$ preserves $(J,V,\omega)$. In particular, $M$ takes the form
$$
M = \pmatrix A\hfill& 0\hfill\cr 0\hfill&B\hfill\cr\endpmatrix,
$$
for some $A,B\in\opSO(2)$. Since $M$ preserves $J_0$, $A=B$, and it follows that $M$ is an element of $\Sigma$, as desired.\qed
\medskip
\noindent The significance of the bilinear form $m$ for the study of convexity is given by the following results.
\proclaim{Lemma \& Definition \nextprocno}
\noindent The restriction of $m$ to any $J$-complex line $L$ in $\Bbb{C}\oplus\Bbb{C}$ is either positive-definite, negative-definite or null. Furthermore, this restriction is null if and only if $L$ intersects $V$ non-trivially, that is, if and only if $\pi_1$ has non-trivial kernel over $L$.
\par
We say that a $J$-complex line $L$ in $\Bbb{C}\oplus\Bbb{C}$ is {\bf positive}, {\bf negative} or {\bf null} according to whether the restriction of $m$ to $L$ is positive-definite, negative-definite or null.
\endproclaim
\proclabel{SecondLexicon}
\proof Let $L\subseteq\Bbb{R}^2\oplus\Bbb{R}^2$ be a $J$-complex line. Since $m(J\cdot,J\cdot)=m$, $m$ has well-defined sign over $L$. Since $L$ is $\omega_k$-lagrangian, for all $\xi,\nu\in L$,
$$
\langle\pi_1(\xi),\pi_2(\nu)\rangle = m(\xi,\nu) + \omega_k(\xi,\nu) = m(\xi,\nu).
$$
It follows that $m$ vanishes over $L$ if and only if one of $\opKer(\pi_1)\minter L$ or $\opKer(\pi_2)\minter L$ is non-trivial. However, since $J_0\pi_1=\pi_2J$, these spaces are either both trivial or both non-trivial, and the result follows.\qed
\proclaim{Lemma \nextprocno, {\bf Convexity\ =\ Positivity}}
\noindent Let $\Omega$ be an open subset of $\Bbb{R}^2$. Let $u:\Omega\rightarrow\Bbb{R}$ be a smooth solution of the Monge--Amp\`ere equation \eqnref{MongeAmpere} and let $S\subseteq\Omega\times\Bbb{R}^2$ denote the graph of its derivative. The function $u$ is strictly convex if and only if every tangent plane of $S$ is positive.
\endproclaim
\proof Indeed, the tangent planes of $S$ are positive if and only if, for all $x\in\Omega$ and for all $\xi\in\Bbb{R}^2$, $\langle\xi,\opHess(u)(x)\xi\rangle > 0$.\qed
\newsubhead{Positivity II - a holomorphic approach}[PositivityIIAHolomorphicApproach]
Define $\phi,\psi:\Bbb{C}\oplus\Bbb{C}\rightarrow\Bbb{C}$ by
$$\eqalign{
\phi(z,w) &:= z+w,\ \text{and}\cr
\psi(z,w) &:= -\overline{z}+\overline{w},\cr}\eqnum{\nexteqnno[SpaceTimeCoordinatesI]}
$$
and define
$$
\Phi := (\phi,\psi).\eqnum{\nexteqnno[SpaceTimeCoordinatesII]}
$$
Note that $\Phi$ is $\Bbb{C}$-linear with respect to the complex structure $J$ of the domain and the complex structure $\hat{J}_0$ of the codomain. In addition, for all $(z,w)\in\Bbb{C}^2$,
$$
(\Phi_*m)((z,w)^t,(z,w)^t) = \frac{1}{2}(\left|z\right|^2 - \left|w\right|^2).\eqnum{\nexteqnno[PushForwardOfMinkowskiMetric]}
$$
We thus consider $\phi$ and $\psi$ as the respective projections of $\Bbb{C}\oplus\Bbb{C}$ onto the space and time components of $m$. This has two useful consequences. The first is the following Lipschitz property.
\proclaim{Lemma \nextprocno}
\noindent Let $P$ be a real plane in $\Bbb{C}\oplus\Bbb{C}$. If
$$
m|_P \geq 0,\eqnum{\nexteqnno[Bilipschitz]}
$$
then the restriction of $\phi$ to $P$ is $\sqrt{2}$-bilipschitz.
\endproclaim
\proclabel{Bilipschitz}
\proof Indeed, for all $(z,w)\in P$,
$$
\left|\phi(z,w)\right| = \left|z+w\right| \leq \sqrt{2}\|z,w\|.
$$
Conversely,
$$
\left|\phi(z,w)\right|^2 \geq \|z,w\|^2 + m((z,w)^t,(z,w)^t) \geq \|z,w\|^2,
$$
as desired.\qed
\medskip
\noindent The second useful consequence is the following conformal description of positivity. First, let $\hat{\Bbb{C}}$ denote the extended complex plane and define the conformal diffeomorphism $\alpha:\Bbb{C}\Bbb{P}^1\rightarrow\hat{\Bbb{C}}$ by
$$
\alpha([z:w]) := \frac{w}{z}.\eqnum{\nexteqnno[ParametrisationOfCP]}
$$
\proclaim{Lemma \nextprocno}
\noindent Let $L$ be a $J$-complex line in $\Bbb{C}\oplus\Bbb{C}$. $L$ is positive if and only if $\left|(\alpha\circ\Phi)(L)\right|<1$ and $L$ is null if and only if $\left|(\alpha\circ\Phi)(L)\right|=1$.
\endproclaim
\proclabel{HolomorphicCharacterisationOfPositiveAndNullPlanes}
\proof Choose $(z,w)\in\Phi(L)$. By \eqnref{PushForwardOfMinkowskiMetric}, $L$ is positive if and only if
$$
\left|w\right|^2 < \left|z\right|^2,
$$
which holds if and only if $\left|\alpha(L)\right|=\left|\alpha([z:w])\right|<1$. Likewise, $L$ is null if and only if
$$
\left|w\right|^2 = \left|z\right|^2,
$$
which holds if and only if $\left|\alpha(L)\right|=1$.\qed
\newsubhead{Bernstein type theorems I}[BernsteinTypeTheorems]
The quaternionic formalism that we have presented here now yields simple proofs of key Bernstein-type theorems. These results are not only interesting in their own right, but also provide the basis for the general compactness theorem that will be developed later.
\proclaim{Theorem \nextprocno}
\noindent If $\Surface$ is a complete, non-negative $J$-holomorphic curve in $\Bbb{C}\oplus\Bbb{C}$, then $\Surface$ is an affine plane.
\endproclaim
\proclabel{JuergensQuaternionicSetting}
\proof Indeed, upon passing to the universal cover, we may assume that $\Surface$ is simply-connected. By Lemma \procref{Bilipschitz}, $\phi$ restricts to a bilipschitz, holomorphic map from $\Surface$ into $\Bbb{C}$. It is thus a conformal diffeomorphism, and $\Surface$ therefore has parabolic conformal type. Now define the meromorphic function $\tau:\Surface\rightarrow\hat{\Bbb{C}}$ by
$$
\tau(x) := (\alpha\circ\Phi)(T_x\Surface).
$$
By Lemma \procref{HolomorphicCharacterisationOfPositiveAndNullPlanes}, $\tau$ takes values in the closed unit disk. By Liouville's theorem, $\tau$ is constant, and $\Surface$ is therefore an affine plane, as desired.\qed
\medskip
\noindent The Bernstein-type theorem \cite{Jorgens} of J\"orgens is an immediate corollary.
\proclaim{Theorem \nextprocno, {\bf J\"orgens}}
\noindent If $u:\Bbb{R}^2\rightarrow\Bbb{R}$ is a smooth solution of the real Monge--Amp\`ere equation \eqnref{MongeAmpere}, then $u$ is a quadratic function.
\endproclaim
\proclabel{JuergensB}
\remark[JurgensB] J\"orgens' theorem was a key result in the study of the real Monge--Amp\`ere equations. It was generalised to Dimensions $3$ and $4$ by Calabi in \cite{Calabi}, and then to arbitrary dimension by Pogorelov in \cite{Pogorelov}. The ideas originated in these papers were later applied to the study of Monge-Amp\`ere equations in diverse settings. It is worth noting, however, that in the higher-dimensional case, there is no useful analogue of the quaternionic approach studied here. For such analogoues, we must instead look to the special lagrangian potential equation, where a rich theory can indeed be developed (see \cite{HarvLawsII} and \cite{SmiSLC}).
\medskip
\proof Indeed, by Lemmas \procref{SolutionsOfMAAreJHolomorphicCurves} and \procref{SecondLexicon}, the graph $\Surface$ of $du$ is a positive $J$-holomorphic curve. In addition, being a graph over $\Bbb{R}^2$, it is complete, and the result now follows by Theorem \procref{JuergensQuaternionicSetting}.\qed
\medskip
Another nice application of Theorem \procref{JuergensQuaternionicSetting} is the following mild improvement of the classification result of complete flat surfaces in $\Bbb{H}^3$ by Volkov--Vladimirova and Sasaki (see Chapter $7.F$ of \cite{Spivak}). We say that an immersed surface $(S,e)$ in $\Bbb{H}^3$ is {\bf flat} whenever it has constant extrinsic curvature equal to $1$.\numberedfootnote{Flat surfaces in $\Bbb{H}^3$ are
studied in terms of their Gauss lifts in \cite{KokRossSajiUmeYam}, \cite{KokRossUmeYam} and \cite{MartMil}. The terminology used in these papers differs from our own: where they write ``weakly complete'', we write ``quasicomplete''; and the objects which they refer to as ``flat fronts'' are for us the Gauss lifts of, possibly singular, flat surfaces.}
\proclaim{Theorem \nextprocno}
\noindent The only quasicomplete, ISC, flat surfaces in $\Bbb{H}^3$ are the level sets of horospheres and the level sets of distance functions to complete geodesics.
\endproclaim
\proclabel{VolkovVladimirovaSasaki}
\remark[VolkovVladimirovaSasaki] When the surface is complete, we recover the result of Volkov--Vladimirova and Sasaki.
\medskip
\proof By Gauss' Theorem, every flat surface is intrinsically flat and is thus locally isometric to $\Bbb{R}^2$. Let $\Omega$ be an open subset of $\Bbb{R}^2$, let $e:\Omega\rightarrow\Bbb{H}^3$ be an isometric immersion, and let $\opII_e$ denote its second fundamental form. By the Codazzi-Mainardi equations, the derivative $(\opD\opII_e)_{ijk}$ of $\opII_e$ is symmetric, so that
$$
\opII_e = \opHess(f),
$$
for some smooth function $f$, unique up to addition of a linear function. We define the immersion $\alpha:\Omega\rightarrow\Bbb{R}^2\oplus\Bbb{R}^2$ by
$$
\alpha(x) := (x,df(x)).
$$
Since $\opDet(\opHess(f))=\opDet(\opII_e)=1$, it follows by Lemma \procref{SolutionsOfMAAreJHolomorphicCurves} that $\alpha$ is $J$-holomorphic. Furthermore, $(\pi_1\circ\alpha)$ is a local isometry and $\alpha^*m=2\opII_e$, so that $\alpha$ is also positive.
\par
Now let $(S,e)$ be a quasicomplete, flat immersed surface in $\Bbb{H}^3$. Upon taking the universal cover, we may suppose that $S$ is simply connected. It then follows by the preceding discussion that there exists a $J$-holomorphic immersion $\alpha:S\rightarrow\Bbb{R}^2\oplus\Bbb{R}^2$ such that $\beta:=(\pi_1\circ\alpha)$ is a local isometry and $\alpha^*m=2\opII_e$. Furthermore, since $e$ is quasicomplete, $\alpha$ is complete, and it follows by Theorem \procref{JuergensQuaternionicSetting} that $(S,\alpha)$ is an affine plane. In particular, $\beta$ defines a global isometry from $S$ into $\Bbb{R}^2$. Upon composing with a rotation if necessary, $e':=e\circ\beta^{-1}$ is an isometric immersion with second fundamental form given by
$$
\opII_{e'} := \pmatrix \lambda\hfill& 0\hfill\cr 0\hfill& 1/\lambda\hfill\cr\endpmatrix,
$$
for some constant $\lambda\in]0,1]$. The result now follows by the fundamental theorem of surface theory: when $\lambda=1$, $S$ is the level set of a horofunction and, when $\lambda<1$, $S$ is the set of points lying at distance $\opArcTanh(\lambda)$ from some complete geodesic in $\Bbb{H}^3$.\qed
\newsubhead{Bernstein type theorems II - non-trivial boundary}[BernsteinTypeTheoremsII]
The case of curves with non-trivial boundary is similar, though more technical.
\proclaim{Theorem \nextprocno}
\noindent Let $P\subseteq\Bbb{C}\oplus\Bbb{C}$ be a plane over which $m$ has signature $(1,0)$. Let $(\Surface,\partial \Surface)\subseteq\Bbb{C}\oplus\Bbb{C}$ be a complete, non-negative $J$-holomorphic curve with boundary. If
\medskip
\myitem{(1)} $\partial S\subseteq P$; and
\medskip
\myitem{(2)} there exists $\epsilon>0$ such that, for all $x\in\partial\Surface$,
$$
m|_{T_xS} \geq \epsilon g|_{T_x\Surface}.\eqnum{\nexteqnno[BoundaryCondition]}
$$
then $\Surface$ is an affine plane.
\endproclaim
\proclabel{BoundaryJuergensQuaternionicSetting}
\proof Upon taking the universal cover, we may suppose that $\Surface$ is simply connected. Let $\Cal{F}$ denote the foliation of $P$ by null lines and let $L\subseteq P$ be any line transverse to this foliation. By Lemma \procref{Bilipschitz}, $\phi$ restricts to a linear isomorphism from $P$ into $\Bbb{C}$. By \eqnref{BoundaryCondition}, there exists $\theta>0$ such that $\phi(\partial\Surface)$ makes an angle of at least $\theta$ with $\phi_*\Cal{F}$ at every point so that, by completeness, $\phi(\partial\Surface)$ is a union of complete Lipschitz graphs over $L$. By Lemma \procref{Bilipschitz} again, the restriction of $\phi$ to $\Surface$ is a conformal diffeomorphism onto its image. In particular $\phi(S)$ is a subset of $\Bbb{C}$ bounded by at most $2$ complete Lipschitz graphs. $(\Surface,\partial\Surface)$ is therefore conformally equivalent to one of $\overline{\Bbb{D}}\setminus\left\{-1\right\}$ or $\overline{\Bbb{D}}\setminus\left\{\pm 1\right\}$, and the Riemann surface $\tilde{\Surface}$, obtained by doubling $\Surface$ along its boundary, is thus conformally equivalent to one of $\Bbb{C}$ or $\Bbb{C}\setminus\left\{0\right\}$. In both cases, $\tilde{\Surface}$ is of parabolic type.
\par
Define the meromorphic function $\tau:S\rightarrow\hat{\Bbb{C}}$ by
$$
\tau(x) := (\alpha\circ\Phi)(T_x\Surface).
$$
By Lemma \procref{HolomorphicCharacterisationOfPositiveAndNullPlanes}, $\tau$ takes values in the closed unit disk. Furthermore, by \eqnref{BoundaryCondition}, there exists $\delta>0$ such that, for all $x\in\partial\Surface$,
$$
\left|\tau(x)\right| \leq 1-\delta.
$$
Since the boundary $\partial\Surface$ of $\Surface$ has at most two isolated singularities, and since $\tau$ is bounded, it follows by the Cauchy integral formula that, for all $x\in\Surface$,
$$
\left|\tau(x)\right| \leq 1-\delta.
$$
\par
For all $r>0$, let $D_r$ and $C_r$ denote respectively the disk and circle in $\Bbb{C}$ of radius $r$ about the origin. By definition of $P$, there exists a circle $C$, which is a strict interior tangent to $C_1$ at some point, such that, for all $x\in\partial\Surface$,
$$
\tau(x)\in C.
$$
Let $R:\hat{\Bbb{C}}\rightarrow\hat{\Bbb{C}}$ denote the conformal reflection through $C$. By the Schwarz reflection principle, $\tau$ extends to a holomorphic function $\tilde{\tau}:\tilde{\Surface}\rightarrow\hat{\Bbb{C}}$ such that, for all $x\in\Surface$,
$$
\tilde{\tau}(\overline{x}) := R\tau(x).
$$
In particular,
$$
\tilde{\tau}(\tilde{S}) \subseteq \overline{D}_{1-\delta}\munion R\overline{D}_{1-\delta}=:\Omega.
$$
Since the complement of $\Omega$ in $\hat{\Bbb{C}}$ has non-trivial interior, it follows by Liouville's theorem that $\tilde{\tau}$, and therefore also $\tau$, is constant. The surface $\Surface$ is thus an affine half-plane, and this completes the proof.\qed
\medskip
\noindent In terms of the real Monge--Amp\`ere equation, this yields the following new Bernstein-type theorem for functions defined over domains with boundary.
\proclaim{Theorem \nextprocno}
\noindent Let $\phi:\Bbb{R}\times[0,\infty[\rightarrow\Bbb{R}$ be a smooth solution of the real Monge--Amp\`ere equation \eqnref{MongeAmpere}. If
\medskip
\myitem{(1)} $\phi$ restricts to a non-trivial quadratic function over $\Bbb{R}\times\left\{0\right\}$; and
\medskip
\myitem{(2)} $\Delta\phi$ is bounded over $\Bbb{R}\times\left\{0\right\}$,
\medskip
\noindent then $\phi$ is a quadratic function.
\endproclaim
\proclabel{BoundaryJuergens}
\proof Indeed, by Lemmas \procref{SolutionsOfMAAreJHolomorphicCurves} and \procref{SecondLexicon}, the graph $\Surface$ of $du$ is a positive $J$-holomorphic curve which satisfies the hypotheses of Theorem \procref{BoundaryJuergensQuaternionicSetting}. The result follows.\qed
\newsubhead{Non-quadratic solutions}[NonQuadraticSolutions]
We conclude this section with the following interesting open problem. It is not clear that Theorems \procref{BoundaryJuergensQuaternionicSetting} and \procref{BoundaryJuergens} are optimal as stated. Although their first conditions are natural from the point of view of Liouville's theorem and the Schwarz reflection principle, their second conditions are not. Nevertheless, the first conditions alone are not sufficient, and the function
$$
\phi(x,y) := \frac{x^2}{y} + \frac{y^3}{12}\eqnum{\nexteqnno[FatalCounterExample]}
$$
shows. Indeed, this function is a convex solution of the real Monge--Amp\`ere equation \eqnref{MongeAmpere} over the open half-plane $\Bbb{R}\times]0,\infty[$ whose restriction to every horizontal line is quadratic, but which is itself trivially not quadratic.
\par
The above function, together with many similar examples, is obtained by directly solving the Monge--Amp\`ere equation \eqnref{MongeAmpere} with the ansatz
$$
f(x,y) = \alpha(y)x^2 + \beta(y)x + \gamma(y).\eqnum{\nexteqnno[Ansatz]}
$$
However, it is also worth reviewing its properties from the quaternionic point of view. Note first that
$$
\Phi(V) = \left\{ (z,\overline{z})\ |\ z\in\Bbb{C}\right\}.\eqnum{\nexteqnno[VerticalSpaceInSpaceTime]}
$$
Let $\Bbb{H}^+$ denote the upper half-space in $\Bbb{C}$ and let $\overline{\Bbb{H}}^+$ denote its closure. Consider a holomorphic function $F:=(f,g)^t:\overline{\Bbb{H}}^+\rightarrow\Bbb{C}\oplus\Bbb{C}$, and suppose that
$$
F(\Bbb{R}) \subseteq \Phi(V).\eqnum{\nexteqnno[CounterExampleBoundaryConditionA]}
$$
By \eqnref{VerticalSpaceInSpaceTime}, for all $x\in\Bbb{R}$,
$$
f(x) = \overline{g(x)}.\eqnum{\nexteqnno[CounterExampleBoundaryConditionB]}
$$
It follows by the Schwarz reflection principle that $f$ and $g$ extend to entire functions such that, for all $z\in\Bbb{H}^+$,
$$
g(z) = \overline{f(\overline{z})}.\eqnum{\nexteqnno[DoublingProperty]}
$$
In particular, it is sufficient to prescribe only the first component $f$. In order to obtain a complete immersed surface with boundary, it is also necessary that $f$ define a conformal diffeomorphism of $\Bbb{H}^+$ onto its image. This holds, for example, when $f$ is quadratic. We thus choose
$$
f(z) := -i-i(z+i)^2,\eqnum{\nexteqnno[FirstComponentOfCounterExample]}
$$
so that
$$
F(z) = (-i-i(z+i)^2,i+i(z-i)^2)^t.\eqnum{\nexteqnno[FormulaForF]}
$$
Composing with $\Phi^{-1}$ yields
$$
(\Phi^{-1}\circ F)(z) = (2xy + 2iy,2x + i(y^2 - x^2))^t,\eqnum{\nexteqnno[TransformedFormulaForF]}
$$
where $x+iy:=z$, so that $(\Phi^{-1}\circ F)$ is the graph over $\Bbb{R}\times]0,\infty[$ of the $1$-form
$$
\alpha(x,y) := \frac{2x}{y}dx + \bigg(\frac{y^2}{4} - \frac{x^2}{y^2}\bigg)dy,\eqnum{\nexteqnno[OneForm]}
$$
which integrates to \eqnref{FatalCounterExample}.
\par
With a view towards the optimal formulations of Theorems \procref{BoundaryJuergensQuaternionicSetting} and \procref{BoundaryJuergens}, it is tempting to conjecture that the only non-quadratic solutions of the Monge--Amp\`ere equation over the half-plane which restrict to quadratic functions over the boundary are those described above. However, although we have not been able to find other examples, we also have no proof of this assertion.
\newhead{Monge--Amp\`ere structures and $k$-surfaces}[MongeAmpereStructuresAndKSurfaces]
\newsubhead{Monge--Amp\`ere structures}[MongeAmpereStructures]
In \cite{LabMA} Labourie uses a non-integrable generalisation of the structures studied in Section \headref{QuaternionsAndBernsteinTypeTheorems} to derive compactness results for families of immersed surfaces of prescribed extrinsic curvature in $3$-dimensional riemannian manifolds. Labourie's work now underlies a large part of the modern theory of such surfaces, with applications in hyperbolic geometry, Teichm\"uller theory, general relativity, and so on. We now provide an elementary presentation of Labourie's ideas.
\par
Let $X$ be a manifold. Labourie defines a {\bf Monge--Amp\`ere structure} over $X$ to be a quadruplet $(W,g,V,J)$ where
\medskip
\myitem{(1)} $W$ is a smooth $4$-dimensional subbundle of the tangent bundle of $X$ furnished with a riemannian metric $g$;
\medskip
\myitem{(2)} $V$ is a smooth, oriented $2$-dimensional subbundle of $W$;
\medskip
\myitem{(3)} $J$ is a smooth section of $\opEnd(W)$ such that, for all $x$, $J_x$ restricts to a compatible complex structure of the fibre $W_x$ with respect to which the fibre $V_x$ is real.
\medskip
\remark[AlternativePerspective] By Lemmas \procref{CharacterisationOfM} and \procref{StabiliserSubgroup} and the subsequent remark, at each point $x$ of $X$, the prescription of $(V_x,J_x)$ is equivalent to the prescription of a compatible quaternionic structure $(I_x,J_x,K_x)$, together with a unit-length symmetric bilinear form $m_x$ satisfying \eqnref{InvarianceOfM}. The above geometric structure can thus be expressed in purely quaternionic terms. This point of view is particularly useful for its extension to the higher-dimensional case (see \cite{SmiSLC}). We also highlight the parallels with the analogous theory of Monge--Amp\`ere structures developed in \cite{LychRub}, \cite{LychRubChek} and \cite{Rub}.
\medskip
The geometric significance of $V$ is easier to understand when described in conformal terms. Indeed, at every point $x$, the set of complex lines in $W_x$ which intersect $V_x$ non-trivially defines a circle $C_x$ in the complex projective space of $(W_x,J_x)$. Furthermore, the orientation of $V_x$ in turn yields an orientation of this circle which therefore has a well-defined interior $D_x$. The prescription of $V$ is thus equivalent to the prescription of a conformal disk bundle $D$ in the complex projective bundle of $(W,J)$. In the spirit of Lemma \procref{HolomorphicCharacterisationOfPositiveAndNullPlanes} this yields a notion of positivity of $J_x$-complex lines in $V_x$. Let $\overline{D}$ denote the bundle whose fibre at $x$ is the topological closure of $D_x$.
\par
A {\bf $J$-holomorphic curve} in $X$ is a pair $(\Surface,\phi)$ where $\Surface$ is a Riemann surface and $\phi:\Surface\rightarrow X$ is a smooth function such that, for all $x$, $\opIm(D\phi(x))\subseteq W_{\phi(x)}$ and
$$
J_{\phi(x)}D\phi(x) = D\phi(x)j_x,\eqnum{\nexteqnno[JHolomorphicCurve]}
$$
where $j$ here denotes the complex structure of $\Surface$. We define a {\bf Monge--Amp\`ere surface} to be an immersed $J$-holomorphic curve $(S,\phi)$ in $X$ such that
\medskip
\myitem{(1)} $(S,\phi)$ is tangent to $W$;
\medskip
\myitem{(2)} $(S,\phi)$ is complete with respect to $g$; and
\medskip
\myitem{(3)} the tangent space to $(S,\phi)$ is at every point an element of $\overline{D}$.
\medskip
\noindent We say that a Monge--Amp\`ere surface is {\bf positive} or {\bf null} at some point according to whether its tangent space at that point is an element of $D$ or $\partial D=\overline{D}\setminus D$. Following Labourie, we define a {\bf curtain surface} to be a Monge--Amp\`ere surface which is everywhere null.
\par
Labourie's framework consists of two parts, namely a compactness result, which we will study in the next section, and a dichotomy result, given by the following theorem.
\proclaim{Theorem \nextprocno}
\noindent Let $(\Surface,\phi)$ be a Monge--Amp\`ere surface and let $(\Surface',\phi')$ be a curtain surface. If there exists $x\in\Surface$ and $x'\in\Surface'$ such that $\phi(x)=\phi'(x')$ and $\opIm(D\phi(x))=\opIm(D\phi'(x'))$, then there exists a neighbourhood $U$ of $x$ in $\Surface$, a neighbourhood $V$ of $x'$ in $\Surface'$, and a conformal diffeomorphism $\alpha:U\rightarrow V$ such that
$$
\phi = \phi'\circ\alpha.\eqnum{\nexteqnno[CurtainSurfaceTheorem]}
$$
\endproclaim
\proclabel{CurtainSurfaceTheorem}
\medskip
\noindent Theorem \procref{CurtainSurfaceTheorem} motivates the following definition. We say that the Monge--Amp\`ere structure is {\bf integrable} whenever every point of $\partial D:=\overline{D}\setminus D$ is tangent to some curtain surface.
\proclaim{Theorem \nextprocno}
\noindent If $(W,V,J)$ is an integrable Monge--Amp\`ere structure then every point of $\overline{D}\setminus D$ is tangent to a unique inextensible, simply-connected curtain surface.
\endproclaim
\proclabel{CurtainSurfaceTheoremI}
\proof This follows immediately from Theorem \procref{CurtainSurfaceTheorem} and unique continuation.\qed
\proclaim{Theorem \nextprocno, {\bf Labourie's Dichotomy}}
\noindent If $(W,V,J)$ is an integrable Monge--Amp\`ere structure, then every Monge--Amp\`ere surface whose tangent bundle meets $\partial D=\overline{D}\setminus D$ at some point is a curtain surface. That is, every Monge--Amp\`ere surface is either everywhere positive, or everywhere null.
\endproclaim
\proclabel{CurtainSurfaceTheoremII}
\proof This likewise follows immediately from Theorem \procref{CurtainSurfaceTheorem} and unique continuation.\qed
\medskip
{\bf\noindent Proof of Theorem \procref{CurtainSurfaceTheorem}:\ }Suppose the contrary. For $r>0$, let $\Bbb{D}_r$ and $\Bbb{R}_r^m$ denote the open disk and the open ball, both of radius $r$, in $\Bbb{C}$ and $\Bbb{R}^m$ respectively. We parametrize a neighbourhood of $y:=\phi(x)=\phi'(x')$ in $X$ by $\Bbb{D}_r\times\Bbb{D}_r\times\Bbb{B}^m_r$ in such a manner that $\Bbb{D}_r\times\left\{(0,0)\right\}$ identifies with a portion of $\phi'(S')$ about $x'$ and, for all $z\in\Bbb{D}_r$,
$$\eqalign{
W_{(z,0,0)} &= \Bbb{C}\times\Bbb{C}\times\left\{0\right\},\ \text{and}\cr
J_{(z,0,0)} &= \pmatrix i\hfill& 0\hfill& 0\hfill\cr 0\hfill& i\hfill& 0\hfill\cr 0\hfill& 0\hfill& 0\hfill\cr\endpmatrix.\cr}
$$
Since $(\Surface,\phi)$ is tangent to $(\Surface',\phi')$ at $y$, a portion of this surface about $x$ identifies with the graph of a smooth function $F:\Bbb{D}_r\rightarrow\Bbb{D}_r\times\Bbb{B}_r^m$ which, by hypothesis, is non-zero. Since both $\Bbb{D}_r\times\left\{(0,0)\right\}$ and the graph of $F$ are $J$-holomorphic curves, it follows by Aronszajn's unique continuation theorem (see \cite{Aron} and Theorem $2.3.4$ of \cite{McDuffSal}) that
$$
F(z) = (P(z),Q(z)) + \opO(z^{k+2}),\eqnum{\nexteqnno[PolynomialApproximation]}
$$
for real homogeneous polynomials $P$ and $Q$, both of order $(k+1)$, where $k\geq 1$. Furthermore, since the graph of $F$ is tangent to $W$,
$$
Q=0,\eqnum{\nexteqnno[QVanishes]}
$$
and, by $J$-holomorphicity again,
$$
P(z) = az^{k+1},\eqnum{\nexteqnno[FormulaForP]}
$$
for some $a\in\Bbb{C}\setminus\left\{0\right\}$.
\par
From this point, the result follows by an elementary topological argument. Let $\Psi$ be a smooth function which, for all $(z,w,\xi)$, sends the complex projective space of $(W_{(z,w,\xi)},J_{(z,w,\xi)})$ to that of $(W_{(z,0,0)},J_{(z,0,0)})$. Since $W$ and $J$ are constant over $\Bbb{D}_t\times\left\{(0,0)\right\}$, we think of $\Psi$ as sending the complex projective bundle of $(W,J)$ to the trivial bundle with fibre $\Bbb{C}\Bbb{P}^1$. For all $z$, let $\phi(z)$ denote the closest point of the boundary of $\Psi(D_{(z,F(z))})$ to the origin. We claim that
$$
\phi(z) = \opO(\left|z\right|^{k+1}).
$$
Indeed, since $\Bbb{D}_r\times\left\{(0,0)\right\}$ is a curtain surface, for all $z\in\Bbb{D}_r$,
$$
0\in\Psi(D_{(z,0,0)}),
$$
and, since $F=\opO(\left|z\right|^{k+1})$, the assertion follows by smoothness of $\Psi$.
\par
Consider now the curves $\gamma_\epsilon(\theta):=\epsilon e^{i\theta}$ and $\delta_\epsilon(\theta):=\phi(\epsilon e^{i\theta})$. For all $\theta$, let $\gamma_\epsilon'(\theta)$ denote the image under $\Psi$ of the tangent plane of the graph of $F$ at the point $(\gamma_\epsilon(\theta),F(\gamma_\epsilon(\theta)))$. By \eqnref{PolynomialApproximation}, \eqnref{QVanishes} and \eqnref{FormulaForP},
$$
\gamma_\epsilon'(\theta) = a(k+1)\epsilon^ke^{ik\theta} + \opO(\epsilon^{k+1}).
$$
It follows that, for sufficiently small $\epsilon$, $(\gamma'_\epsilon-\delta_\epsilon)$ winds $k$ times around the origin. However, this is absurd, since the tangent plane of $(S,\phi)$ is at every point an element of $D$ whilst $\delta_\epsilon$ lies along its boundary. The result follows.\qed
\newsubhead{Compactness}[Compactness]
We now prove a compactness property for families of Monge--Amp\`ere surfaces. In order to correctly express the result, we require the concept of smooth Cheeger--Gromov convergence. First, we define a {\bf pointed riemannian manifold} to be a triplet $(X,g,p)$, where $X$ is a smooth manifold, $g$ is a riemannian metric and $p$ is a point of $X$. We say that a sequence $(X_m,g_m,p_m)$ of complete, pointed riemannian manifolds converges to the complete, pointed riemannian manifold $(X_\infty,g_\infty,p_\infty)$ in the {\bf smooth Cheeger--Gromov sense} whenever there exists a sequence $(\Phi_m)$ of functions such that
\medskip
\myitem{(1)} for all $m$, $\Phi_m:X_\infty\rightarrow X_m$ and $\Phi_\infty(p_\infty)=p_m$; and
\medskip
\noindent for every relatively compact open subset $\Omega$ of $X_\infty$, there exists $M$ such that
\medskip
\myitem{(2)} for all $m\geq M$, the restriction of $\Phi_m$ to $\Omega$ defines a smooth diffeomorphism onto its image; and
\medskip
\myitem{(3)} the sequence $((\Phi_m|_\Omega)^*g_m)_{m\geq M}$ converges to $g_\infty|_\Omega$ in the $C^\infty$ sense.
\medskip
\noindent We call $(\Phi_m)$ a sequence of {\bf convergence maps} of $(X_m,g_m,p_m)$ with respect to $(X_\infty,g_\infty,p_\infty)$.
\par
It is reassuring to observe that it defines a Hausdorff topology over the set of isometry classes of complete, pointed riemannian manifolds.\numberedfootnote{Strictly speaking, the class of all complete, pointed riemannian manifolds is not a set. However, the class of all embedded $n$-dimensional submanifolds of $\Bbb{R}^{2n+2}$ furnished with complete riemannian metrics is. Since Whitney's theorem identifies these two classes up to isometry equivalence, the family of isometry equivalence classes of complete, pointed riemannian manifolds can be legitimately considered a set.} Furthermore, although the convergence maps are trivially non-unique, any two sequences $(\Phi_m)$ and $(\Phi_m')$ of convergence maps are equivalent in the sense that there exists an isometry $\Psi:X_\infty\rightarrow X_\infty$ preserving $p_\infty$ such that, for any two relatively compact open subsets $U\subseteq\overline{U}\subseteq V$ of $X_\infty$, there exists $M$ such that
\medskip
\myitem{(1)} for all $m\geq M$, the respective restrictions of $\Phi_m$ and $\Phi_m'\circ\Psi$ to $U$ and $V$ define smooth diffeomorphisms onto their images;
\medskip
\myitem{(2)} for all $m\geq M$, $(\Phi_m'\circ\Psi)(U)\subseteq\Phi_m(V)$; and
\medskip
\myitem{(3)} the sequence $((\Phi_m|_V)^{-1}\circ\Phi_m'\circ\Psi)$ converges in the $C^\infty$ sense to the identity map over $U$.
\medskip
\noindent This condition trivially defines an equivalence relation over the space of sequences of convergence maps.
\par
We now return to the case of immersed submanifolds. We say that a sequence $(\Surface_m,e_m,p_m)$ of complete pointed immersed submanifolds in a complete riemannian manifold $(X,g)$ converges to the complete pointed immersed submanifold $(\Surface_\infty,e_\infty,p_\infty)$ in the smooth Cheeger--Gromov sense whenever the sequence $(\Surface_m,e_m^*g,p_m)$ converges to $(\Surface_\infty,e_\infty^*g,p_\infty)$ in the smooth Cheeger--Gromov sense and, for one, and therefore for any, sequence $(\Phi_m)$ of convergence maps, the sequence $(e_m\circ\Phi_m)$ converges to $e_\infty$ in the $C^\infty_\oploc$ sense.
\par
Smooth Cheeger--Gromov convergence can also be characterised in terms of graphs. Indeed, let $N\Surface_\infty$ denote the normal bundle of $(\Surface_\infty,e_\infty)$ in $\phi_\infty^*TX$. Recall that the exponential map of $X$ defines a smooth function $\opExp:N\Surface_\infty\rightarrow X$. In particular, given a sufficently small smooth section $f:\Omega\rightarrow N\Surface_\infty$ defined over an open subset $\Omega$ of $\Surface_\infty$, the composition $\opExp\circ f$ defines a smooth immersion of $\Omega$ in $X$ which we call the {\bf graph} of $f$. A straightforward but technical argument shows that the sequence $(\Surface_m,e_m,p_m)$ converges to $(\Surface_\infty,e_\infty,p_\infty)$ in the smooth Cheeger--Gromov sense if and only if there exists a sequence $(p_m')$ of points in $\Surface_\infty$ and sequences of functions $(f_m)$ and $(\alpha_m)$ such that
\medskip
\myitem{(1)} $(p_m')$ converges to $p_\infty$;
\medskip
\myitem{(2)} for all $m$, $f_m$ maps $\Surface_\infty$ into $N\Surface_\infty$, $\alpha$ maps $\Surface_\infty$ into $\Surface_m$ and $\alpha(p_m')=p_m$; and
\medskip
\noindent for every relatively compact open subset $\Omega$ of $\Surface_\infty$, there exists $M$ such that
\medskip
\myitem{(3)} for all $m\geq M$, the restriction of $f_m$ to $\Omega$ defines a smooth section of $N\Surface_\infty$ over this set, the restriction of $\alpha_m$ to $\Omega$ defines a smooth diffeomorphism onto its image, and
$$
\opExp\circ f_m|_\Omega = e_m\circ\alpha|_\Omega;\ \text{and}
$$
\myitem{(4)} the sequence $(f_m)_{m\geq M}$ tends to zero in the $C^\infty$ sense.
\medskip
These definitions are readily extended in a number of ways. For example, in the case of a sequence $(X_m,g_m,p_m)$ of pointed riemannian manifolds, the hypothesis of completeness is unnecessary. Instead, it is sufficient to assume that for all $R>0$, there exists $M$ such that, for all $m\geq M$, the closed ball of radius $R$ about $p_m$ in $(X_m,g_m)$ is compact. Likewise, in the case of immersed submanifolds, the target space can be replaced with a sequence $(X_m,g_m,p_m)$ of pointed riemannian manifolds converging in the smooth Cheeger--Gromov sense to some complete pointed riemannian manifold. Furthermore, it is not necessary to suppose that the riemannian manifolds in this sequence are complete, and so on.
\par
Elliptic regularity, together with the Arzela-Ascoli Theorem of \cite{SmiAA}, yields the following compactness result.
\proclaim{Lemma \nextprocno}
\noindent Let $(X,g)$ be a complete riemannian manifold furnished with a Monge--Amp\`ere structure $(W,V,J)$. Let $(\Surface_m,e_m,p_m)$ be a sequence of complete pointed $J$-holomorphic curves in $(X,g)$. If
\medskip
\myitem{(1)} there exists a compact subset $K$ of $X$ such that $e_m(p_m)\in X$ for all $m$; and
\medskip
\myitem{(2)} for every compact subset $L\subseteq X$, there exists a constant $B>0$ such that, for all $m$, the norm of the second fundamental form of $e_m$ is at every point of $e_m^{-1}(L)$ bounded above by $B$,
\medskip
\noindent then there exists a complete pointed Monge--Amp\`ere surface $(\Surface_\infty,e_\infty,p_\infty)$ in $(X,g)$ towards which this sequence subconverges in the smooth Cheeger--Gromov sense.
\endproclaim
\proclabel{Precompactness}
\remark[Precompactness] This result readily extends in many ways. For example, it is not necessary to assume that the $J$-holomorphic curves are complete. Likewise, the target manifold with its Monge--Amp\`ere structure may be replaced by sequences converging to some limit, and so on.
\medskip
\noindent We now recall a simplified version of the {\bf quasi-maximum lemma} (see \cite{Gromov}).
\proclaim{Lemma \nextprocno}
\noindent Let $X$ be a complete metric space and let $f:X\rightarrow[1,\infty[$ be an upper semi-continuous function. For all $x\in X$, there exists $y\in X$ such that
\medskip
\myitem{(1)} $f(y)\geqslant f(x)$; and
\medskip
\myitem{(2)} for all $z\in B_{1/\sqrt{f(y)}}(y)$, $f(z)\leqslant 2f(y)$.
\endproclaim
\proclabel{QuasiMaximumLemma}
\proof Indeed, otherwise, we recursively construct a sequence $(x_m)$ of points of $X$ such that, for all $m$, $f(x_{m+1})\geqslant 2f(x_m)$ and $d(x_{m+1},x_m)\leqslant 1/\sqrt{f(x_m)}$. We readily verify that $(x_m)$ is a Cauchy sequence which therefore converges to some limit $x_\infty$, say. However, by upper semi-continuity,
$$
f(x_\infty) \geqslant \mlimsup_{m\rightarrow\infty} f(x_m) = \infty,
$$
which is absurd, and the result follows.\qed
\proclaim{Lemma \nextprocno}
\noindent Let $(X,g)$ be a complete riemannian manifold furnished with a Monge--Amp\`ere structure $(W,V,J)$. For every compact subset $K$ of $X$, there exists $B>0$ such that every Monge--Amp\`ere surface $(\Surface,e)$ in $X$ has second fundamental form bounded above by $B$ at every point of $e^{-1}(K)$.
\endproclaim
\proclabel{Compactness}
\proof We prove this result using a standard blow-up argument. Indeed, suppose the contrary. There exists a sequence $(\Surface_m,e_m,p_m)$ of complete, pointed Monge--Amp\`ere surfaces and a sequence $(B_m)$ of positive numbers converging to $+\infty$ such that, for all $m$, $e_m(p_m)\in K$ and
$$
\|\opII_m(p_m)\| = B_m,
$$
where, for all $m$, $\opII_m$ denotes the second fundamental form of $e_m$. By Lemma \procref{QuasiMaximumLemma}, we may suppose that, for all $m$ and for all $y\in B_{1/\sqrt{B_m}}(p_m)$,
$$
\|\opII_m(y)\|\leq 2B_m.
$$
We extend $g$ to a riemannian metric defined over the whole of $X$ and, upon rescaling, we consider the sequence $(X,B_m^2g,e_m(p_m))$ of pointed riemannian manifolds. This sequence converges in the Cheeger--Gromov sense to $(X_\infty,g_\infty,0)$, where $X_\infty:=\Bbb{C}\times\Bbb{C}\times\Bbb{R}^m$ and $g_\infty$ is the standard euclidean metric. In addition, we may suppose that the Monge--Amp\`ere structures of these manifolds converge to the constant Monge--Amp\`ere structure $(W_\infty,V_\infty,J_\infty)$ given by
$$\eqalign{
W_\infty&:=\Bbb{C}\times\Bbb{C}\times\left\{0\right\},\vphantom{\bigg)}\cr
V_\infty&:=\Bbb{R}\times\Bbb{R}\times\left\{0\right\},\ \text{and}\vphantom{\bigg)}\cr
J_\infty&:=\pmatrix 0\hfill& i\hfill\cr i\hfill& 0\hfill\cr\endpmatrix.\cr}
$$
For all $m$, let $\opII_m'$ denote the shape operator of $e_m$ with respect to the rescaled metric $B_m^2g$. For all $m$,
$$
\|\opII_m'(p_m)\| = 1,
$$
and, for all $y\in B_{\sqrt{B_m}}(x_m)$, $\|\opII_m'(y)\|\leq 2$. It follows by Theorem \procref{Precompactness} and the subsequent remark that there exists a Monge--Amp\`ere surface $(\Surface_\infty,e_\infty,p_\infty)$ in $X_\infty$ towards which $(\Surface_m,e_m,p_m)$ subconverges in the smooth Cheeger--Gromov sense. In particular, the norm of the second fundamental form of this surface at $p_\infty$ is equal to $1$. However, since this surface is contained in $\Bbb{C}\times\Bbb{C}\times\left\{0\right\}$, by Theorem \procref{JuergensQuaternionicSetting}, it is an affine plane. This yields the desired contradiction and the result follows.\qed
\medskip
\noindent Combining Lemmas \procref{Precompactness} and \procref{Compactness} yields Labourie's compactness theorem.
\proclaim{Theorem \nextprocno, {\bf Labourie's compactness theorem}}
\noindent Let $(X,g)$ be a complete riemannian manifold furnished with a Monge--Amp\`ere structure $(W,V,J)$. Let $(S_m,e_m,p_m)_\minn$ be a sequence of pointed Monge--Amp\`ere surfaces in $(X,g)$. If $(e_m(p_m))_{\minn}$ is precompact in $X$, then $(S_m,e_m,p_m)_\minn$ is precompact in the smooth Cheeger--Gromov topology.
\endproclaim
\proclabel{LabouriesCompactnessTheorem}
\newsubhead{Applications I - $k$-surfaces}[ApplicationsIKSurfaces]
Let $(X,h)$ be a complete, oriented $3$-dimensional riemannian manifold, let $SX$ denote its bundle of unit tangent spheres, and furnish this manifold with the Sasaki metric. We construct an integrable Monge--Amp\`ere structure over $SX$ as follows. First, let $VSX$ denote the vertical subbundle of $TSX$, and let $HSX$ denote the horizontal bundle of the Levi--Civita connection. Recall that, at the point $\xi:=\xi_p\in S X$, we have the natural identifications
$$\eqalign{
H_\xi SX &= T_pX,\ \text{and}\cr
V_\xi SX &= \langle\xi\rangle^\perp,\cr}\eqnum{\nexteqnno[HorizontalAndVerticalSubbundles]}
$$
where $\langle\xi\rangle^\perp$ here denotes the orthogonal complement of $\xi_p$ in $T_pX$, so that $T_\xi SX$ decomposes as
$$
T_\xi SX = T_pX \oplus \langle\xi_p\rangle^\perp.\eqnum{\nexteqnno[DecompositionOfTangentBundle]}
$$
With respect to this decomposition, we define the subbundles $W$ and $V$ by
$$\eqalign{
W_\xi &:= \langle\xi_p\rangle^\perp\oplus\langle\xi_p\rangle^\perp,\ \text{and}\cr
V_\xi &:= \left\{0\right\}\oplus\langle\xi_p\rangle^\perp.\cr}\eqnum{\nexteqnno[MAWAndV]}
$$
Since $X$ is oriented, its tangent bundle carries a well-defined wedge product, given by
$$
h(\nu_p,\nu_p'\wedge\nu_p'') := \opdVol(\nu_p,\nu_p',\nu_p'').\eqnum{\nexteqnno[DefinitionOfWedgeProduct]}
$$
Thus, for all $\xi:=\xi_p\in SX$, $\langle\xi\rangle^\perp$ carries a well-defined complex structure $j_\xi$ given by
$$
j_\xi\cdot\nu_p = \xi_p\wedge\nu_p.\eqnum{\nexteqnno[ComplexStructureOverV]}
$$
Given a smooth function $\phi:SX\rightarrow\Bbb{R}$, we define the complex structure $J_\phi$ over $W$ by
$$
J_{\phi,\xi}\cdot(\nu_p,\mu_p) := (e^\phi j_\xi\mu_p,e^{-\phi} j_\xi\nu_p).\eqnum{\nexteqnno[MAComplexStructure]}
$$
Since, in addition, $j_\xi$ defines an orientation over $V$, $(W,V,J_\phi)$ defines a Monge--Amp\`ere structure over $X$. Note that the symmetric bilinear form $m$ corresponding to this structure is, up to sign, given by
$$
m_\xi((\nu_p,\mu_p),(\nu_p,\mu_p)) = 2\langle\nu_p,\mu_p\rangle.\eqnum{\nexteqnno[MApplicationI]}
$$
\par
We now verify integrability of this Monge--Amp\`ere structure. For every complete geodesic $\Gamma$ in $X$, let $\opN\Gamma$ denote its unit normal bundle in $SX$.
\proclaim{Lemma \nextprocno}
\noindent The Monge--Amp\`ere structure $(W,V,J_\phi)$ is integrable and its curtain surfaces are the covers of those surfaces of the form $\opN\Gamma$, for some complete geodesic $\Gamma$ in $X$.
\endproclaim
\proclabel{IntegrabilityI}
\proof Let $\Gamma$ be a complete geodesic in $X$. With respect to the decomposition \eqnref{DecompositionOfTangentBundle}, given any $\xi:=\xi_p\in\opN_p\Gamma$, the tangent space to $\opN\Gamma$ at this point is given by
$$
\opT_\xi\opN\Gamma = \langle (\tau,0),(0,\tau\wedge\xi)\rangle,\eqnum{\nexteqnno[TangentSpaceOfNGamma]}
$$
where $\tau$ is here a unit tangent vector to $\Gamma$ at $p$. It immediately follows that $\opN\Gamma$ is a curtain surface. Conversely, given a point $\xi:=\xi_p$ in $SX$ and a vertical vector $(0,\nu_p)\in W_\xi$, let $\Gamma$ denote the unique geodesic tangent to $\xi_p\wedge\nu_p$ at $p$. By \eqnref{TangentSpaceOfNGamma}, $\opN\Gamma$ is a curtain surface tangent to $(0,\nu_p)$ at $\xi_p$. Since $\xi$ is arbitrary, integrability follows, and this completes the proof.\qed
\medskip
It remains only to study the geometry of positive Monge--Amp\`ere surfaces. First, let $(\Surface,e)$ be an ISC immersed surface in $X$ and let $(\Surface,\hat{e})$ denote its Gauss lift. Recall from the introduction that, for any smooth function $\kappa:SX\rightarrow\Bbb{R}$, the extrinsic curvature $K_e$ of $(\Surface,e)$ is said to be {\bf prescribed} by $\kappa$ whenever
$$
K_e = \kappa\circ\hat{e}.\eqnum{\nexteqnno[PrescribedCurvature]}
$$
\proclaim{Lemma \nextprocno}
\noindent An immersed surface in $SX$ is a positive Monge--Amp\`ere surface if and only if it is the Gauss lift of a quasicomplete, ISC, immersed surface in $X$ of extrinsic curvature prescribed by $e^{2\phi}$.
\endproclaim
\proclabel{DescriptionOfMongeAmpereSurfacesI}
\proof We first show that every positive Monge--Amp\`ere surface is the Gauss lift of some immersed surface. Indeed, since $V$ is the kernel of the canonical projection $\pi:SX\rightarrow X$, for every positive Monge--Amp\`ere surface $(\Surface,\hat{e})$, the function $e:=\pi\circ\hat{e}$ is an immersion. Observe, furthermore, that $\hat{e}$ is normal to $e$, so that the orientation of $\Surface$ can be chosen in such a manner that $(\Surface,\hat{e})$ is the Gauss lift of $(S,e)$, as desired.
\par
Now let $(\Surface,e)$ be an oriented, immersed surface in $X$ and let $(\Surface,\hat{e})$ denote its Gauss lift. Trivially $(S,e)$ is quasicomplete if and only if $(\Surface,\hat{e})$ is complete. Next, with respect to the decomposition \eqnref{DecompositionOfTangentBundle}, the derivative of $\hat{e}$ at a point $p\in\Surface$ is given by
$$
D\hat{e}(p)\cdot\xi_p = (De(p)\cdot\xi_p,De(p)\cdot A_e(p)\cdot \xi_p).
$$
Thus, by Lemma \procref{FirstLexicon}, $(S,\hat{e})$ is a $J$-holomorphic curve if and only if $(\Surface,e)$ has extrinsic curvature prescribed by $e^{2\phi}$. Likewise, by \eqnref{MApplicationI}, $(S,\hat{e})$ is positive if and only if $A_e$ is everywhere positive-definite, that is, if and only if $(S,e)$ is ISC. The result now follows.\qed
\medskip
\noindent The theory developed in Sections \subheadref{MongeAmpereStructures} and \subheadref{Compactness} now yields the following result.
\proclaim{Theorem \nextprocno, {\bf Labourie}}
\noindent Let $\phi:SX\rightarrow\Bbb{R}$ be a smooth function. Let $(S_m,e_m,p_m)$ be a sequence of quasicomplete, pointed, immersed surfaces of extrinsic curvature prescribed by $e^{2\phi}$. If $(e_m(p_m))$ is precompact in $X$, then the sequence $(S_m,\hat{e}_m,p_m)$ of Gauss lifts is precompact in the smooth Cheeger--Gromov sense. Furthermore, any accumulation point $(S_\infty,\hat{e}_\infty,p_\infty)$ of this sequence which is not a curtain surface is the Gauss lift of some quasicomplete, pointed immersed surface of constant extrinsic curvature prescribed by $e^{2\phi}$.
\endproclaim
\proclabel{LabouriesCompactnessTheoremPGCI}
\remark[LabouriesCompactnessTheoremPGCI] As in Section \subheadref{Compactness}, we have stated this result in its simplest form. More generally, the ambient space $X$ and that function $\phi$ may be allowed to vary, the surfaces need not be quasicomplete, and so on.
\newsubhead{Applications II - Isometric immersions}[ApplicationsIIIsometricImmersions]
We conclude this paper by reviewing a more sophisticated version of the construction of the preceding section developed in \cite{LabMP}. By incorporating the domain explicitly into the construction, we are able to study families of isometric immersions, as opposed to reparametrisation equivalence classes, as is the case in the preceding section.
\par
We first review the corresponding linear framework. Let $\opIsom(\Bbb{R}^2,\Bbb{R}^3)$ denote the subset of $\opLin(\Bbb{R}^2,\Bbb{R}^3)$ consisting of linear isometries $\alpha:\Bbb{R}^2\rightarrow\Bbb{R}^3$. Note that this is an orbit of $\opSO(3)$ in $\opLin(\Bbb{R}^2,\Bbb{R}^3)$ upon which this group acts with trivial stabiliser. It is therefore a smooth $3$-dimensional submanifold, and its tangent bundle admits trivialisation
$$
\tau:\opLin(\Bbb{R}^2,\Bbb{R}^3)\times\opso(3)\rightarrow\opT\opIsom(\Bbb{R}^2,\Bbb{R}^3);(\alpha,A)\mapsto(\alpha,A\circ\alpha).\eqnum{\nexteqnno[TrivialisationI]}
$$
In addition, viewing the wedge product linear isomorphism from $\Bbb{R}^3$ to $\opso(3)$, a useful alternative form of this trivialisation is
$$
\sigma:\opLin(\Bbb{R}^2,\Bbb{R}^3)\times\Bbb{R}^3\rightarrow\opT\opIsom(\Bbb{R}^2,\Bbb{R}^3);(\alpha,x)\mapsto(\alpha,x\wedge\alpha).\eqnum{\nexteqnno[TrivialisationII]}
$$
\par
We now return to the non-linear case. Let $(\Surface,g)$ be an oriented riemannian surface, and let $(X,h)$ be an oriented $3$-dimensional riemannian manifold. Over the cartesian product $S\times X$, we consider the bundle $\opLin(\pi_1^*T\Surface,\pi_2^*TX)$, where $\pi_1$ and $\pi_2$ respectively denote the projections onto the first and second factors. Let $\opIsom(S\times X)$ denote its subbundle whose fibre at the point $(p,q)$ is $\opIsom(T_p\Surface,T_qX)$. We construct a Monge--Amp\`ere structure over the total space of this subbundle as follows. Let $\opV\opIsom(\Surface\times X)$ denote the vertical subbundle of $\opT\opIsom(\Surface\times X)$, and let $\opH\opIsom(\Surface\times X)$ denote the horizontal bundle of the Levi-Civita connection. Bearing in mind the discussion of the preceding paragraph, at the point $\alpha:=\alpha_{(p,q)}\in\opIsom_{(p,q)}(\Surface\times X)$, we have the natural identifications,
$$\eqalign{
H_\alpha\opIsom(\Surface\times X) &:= T_pX\oplus T_qX,\ \text{and}\cr
V_\alpha\opIsom(\Surface\times X) &:= T_qX,\cr}\eqnum{\nexteqnno[HorizontalAndVerticalSubbundlesII]}
$$
so that $T_\alpha\opIsom(\Surface\times X)$ decomposes as
$$
\opT_\alpha\opIsom(\Surface\times X) = T_p\Surface\oplus T_qX\oplus T_q X.\eqnum{\nexteqnno[DecompositionOfTIsom]}
$$
With respect to this decomposition, we define the subbundles $W$ and $V$ by
$$\eqalign{
W_\alpha &:= \left\{ (\xi_p,\alpha\cdot\xi_p,\alpha\cdot\nu_p)\ |\ \xi_p,\nu_p\in T_p\Surface\right\},\ \text{and}\cr
V_\alpha &:= \left\{ (0,0,\alpha\cdot\nu_p)\ |\ \nu_p\in T_p\Surface\right\}.\cr}\eqnum{\nexteqnno[DefinitionOfWAndVII]}
$$
We furnish $W$ with the Sasaki metric $g$, and, given a smooth function $\phi:\opIsom(\Surface\times X)\rightarrow\Bbb{R}$, we define the complex structure $J_\phi$ over $W$ by
$$
J_{\phi,\alpha}\cdot(\xi_p,\alpha\cdot\xi_p,\alpha\cdot\nu_p) := (e^\phi j_p\cdot\nu_p,e^\phi \alpha\cdot j_p\cdot\nu_p,e^{-\phi}\alpha\cdot j_p\cdot\xi_p),\eqnum{\nexteqnno[DefinitionOfJII]}
$$
where $j$ here denotes the complex structure of $(\Surface,g)$. Since the orientation of $\Surface$ yields an orientation of $V$, the quadruplet $(W,g,V,J_\phi)$ again defines a Monge--Amp\`ere structure over $\opIsom(\Surface\times X)$. As before, we verify that the associated symmetric, bilinear form $m$ is given by
$$
m_\alpha((\xi_p,\alpha\cdot\xi_p,\alpha\cdot\nu_p),(\xi_p,\alpha\cdot\xi_p,\alpha\cdot\nu_p)) := 2\langle\xi_p,\nu_p\rangle.\eqnum{\nexteqnno[DefinitionMII]}
$$
\par
We now show that this Monge--Amp\`ere structure is integrable. To this end, let $\gamma:\Bbb{R}\rightarrow\Surface$ and $\delta:\Bbb{R}\rightarrow X$ be unit speed parametrised geodesics and let $\alpha:T_{\gamma(0)}\Surface\rightarrow T_{\delta(0)}X$ be an isometry such that
$$
\alpha\cdot\dot{\gamma}(0) = \dot{\delta}(0).\eqnum{\nexteqnno[PropertiesOfAlpha]}
$$
We define $e_{\gamma,\delta,\alpha}:\Bbb{R}^2\rightarrow\opIsom(\Surface\times X)$ by
$$
e_{\gamma,\delta,\alpha}(s,t) := (\gamma(s),\delta(s),\opExp(t\dot{\delta}(s)\wedge)\cdot\tau_s\alpha),\eqnum{\nexteqnno[DefinitionOfCurtainSurfaceII]}
$$
where, for all $s$, $\tau_s$ here denotes parallel transport along $(\gamma,\delta)$ from $(\gamma(0),\delta(0))$ to $(\gamma(s),\delta(s))$.
\proclaim{Lemma \nextprocno}
\noindent $(W,V,J_\phi)$ is integrable, and its curtain surfaces are quotients of immersed surfaces of the form $(\Bbb{R}^2,e_{\gamma,\delta,\alpha})$, with $\gamma$, $\delta$ and $\alpha$ as above.
\endproclaim
\proof Indeed, with respect to the decomposition \eqnref{DecompositionOfTIsom},
$$\eqalign{
De_{\gamma,\delta,\alpha}(s,t)\cdot\partial_s &= (\dot{\gamma}(s),\alpha\cdot\dot{\gamma}(s),0),\ \text{and}\cr
De_{\gamma,\delta,\alpha}(s,t)\cdot\partial_t &= (0,0,\alpha\cdot\dot{\gamma}(s)),\cr}\eqnum{\nexteqnno[DerivativeOfCurtainII]}
$$
so that $(\Bbb{R}^2,e_{\gamma,\delta,\alpha})$ is a curtain surface. Conversely, given a point $\alpha:=\alpha_{(p,q)}\in\opIsom(\Surface\times X)$ and a vertical vector $(0,0,\alpha\cdot\xi_p)$ at this point, let $\gamma:\Bbb{R}\rightarrow\Surface$ and $\delta:\Bbb{R}\rightarrow X$ denote the unique unit speed geodesics such that
$$\eqalign{
\dot{\gamma}(0) &:= \xi_p,\ \text{and}\cr
\dot{\delta}(0) &:= \alpha\cdot\xi_p.\cr}
$$
By \eqnref{DerivativeOfCurtainII}, $(\Bbb{R}^2,e_{\gamma,\delta,\alpha})$ is a curtain surface tangent to $(0,0,\alpha\xi_p)$ at $\alpha_{(p,q)}$. Integrability follows, and this completes the proof.\qed
\medskip
It remains only to describe the geometric properties of positive Monge--Amp\`ere surfaces. Let $(\Surface',g')$ be another riemannian surface and consider an immersion $(e,f):\Surface'\rightarrow\Surface\times X$. We say that $(\Surface',(e,f))$ is {\bf admissable} whenever both $e$ and $f$ are isometries. When this holds, we define $E_{e,f}:\Surface'\rightarrow\opIsom(\Surface\times X)$ by
$$
E_{e,f} := (e,f,Df\circ De^{-1}).\eqnum{\nexteqnno[GaussLiftII]}
$$
We call $(\Surface',E_{e,f})$ the {\bf Gauss lift} of $(\Surface',(e,f))$. We say that $(\Surface',(e,f))$ is {\bf infinitesimally strictly complex} and {\bf quasicomplete} whenever $(\Surface',f)$ has these properties and, given a smooth function $\kappa:\opIsom(\Surface\times X)\rightarrow\Bbb{R}$, we say that $(\Surface',(e,f))$ has extrinsic curvature {\bf prescribed} by $\kappa$ whenever the extrinsic curvature $K_f$ of $f$ satisfies
$$
K_f = \kappa\circ E_{e,f}.\eqnum{\nexteqnno[PrescribedCurvatureII]}
$$
\proclaim{Lemma \nextprocno}
\noindent An immersed surface in $\opIsom(\Surface\times X)$ is a positive Monge--Amp\`ere surface if and only if it is the Gauss lift of a quasicomplete, ISC, admissable immersed surface $(\Surface',(e,f))$ in $\Surface\times X$ of extrinsic curvature prescribed by $e^{2\phi}$.
\endproclaim
\proclabel{DescriptionOfMongeAmpereSurfacesII}
\proof As in the proof of Lemma \procref{DescriptionOfMongeAmpereSurfacesI}, every positive Monge--Amp\`ere surface is the Gauss lift of some immersed surface. Now let $(\Surface',(e,f))$ be an admissable immersed surface and let $(\Surface',E_{e,f})$ denote its Gauss lift. Let $A_f$ denote the shape operator of $f$. Bearing in mind \eqnref{GaussLiftII}, with respect to the decomposition \eqnref{DecompositionOfTIsom}, the derivative of $E_{e,f}$ at a point $p\in\Surface'$ is given by
$$\eqalign{
DE_{e,f}(p)\cdot\xi_p &= (De(p)\cdot\xi_p,Df(p)\cdot\xi_p,Df(p)\cdot A_f(p)\cdot\xi_p)\cr
&= (De(p)\cdot\xi_p,\alpha\cdot De(p)\cdot\xi_p,\alpha\cdot De(p)\cdot A_f(p)\cdot\xi_p),\cr}
$$
where $\alpha:=E_{e,f}(p)$. It follows by Lemma \procref{FirstLexicon} that $(\Surface',E_{e,f})$ is a $J$-holomorphic curve if and only if $(\Surface',(e,f))$ has extrinsic curvature prescribed by $e^{2\phi}$. Likewise, by \eqnref{DefinitionMII}, $(S,\hat{e})$ is positive if and only if $A_e$ is everywhere positive-definite, that is, if and only if $(S',(e,f))$ is ISC. The result now follows.\qed
\medskip
This construction can be applied to the study of isometric immersions in a variety of ways. Since the general result is technical and unenlightening, we shall content ourselves here with an illustrative example. Let $M:=(M,h)$ be a $3$-dimensional Cartan-Hadamard manifold, and let $\Bbb{H}^2$ denote $2$-dimensional hyperbolic space with its metric $g$. Let $(\omega_m)$ be a sequence of bounded, smooth functions over $\Bbb{H}^2$ converging in the $C^\infty_\oploc$ sense to the bounded, smooth function $\omega_\infty$. For all $m$, denote $g_m:=e^{2\omega_m}g$, $S_m:=(\Bbb{H}^2,g_m)$, $X_m:=\opIsom(S_m,M)$, and let $\phi_m:X_m\rightarrow\Bbb{R}$ be such that, for all $\alpha:=\alpha_{p,q}$,
$$
e^{2\phi_m(\alpha)} = \kappa_m(p) - \sigma(\opIm(\alpha)),\eqnum{\nexteqnno[DefinitionOfPhiFinalExample]}
$$
where $\kappa_m$ here denotes the curvature of $g_m$, and $\sigma$ denotes the sectional curvature of $h$. For all $m\in\Bbb{N}$, let $e_m:S_m\rightarrow M$ be an isometric immersion and denote
$$
\hat{e}_m(x) := (x,e_m(x)).\eqnum{\nexteqnno[DefinitionOfHatEm]}
$$
Note that, for all $m$, $\hat{e}_m$ is admissable, quasicomplete, ISC, and of extrinsic curvature prescribed by $e^{2\omega_m}$. Thus, for all $m$, the Gauss lift $E_m$ of $\hat{E}_m$ is a $J_{\phi_m}$-holomorphic curve. The theory developed in Sections \subheadref{MongeAmpereStructures} and \subheadref{Compactness} now yields the following result.
\proclaim{Theorem \nextprocno, {\bf Labourie}}
\noindent Let $x_0\in\Bbb{H}^2$ be a fixed point. If $(e_m(x_0))$ is precompact, then the sequence $(S_m,E_m,x_0)$ of Gauss lifts is precompact. Furthermore, any accumulation point $(S_\infty,E_\infty,x_0)$ is either a tube, or the Gauss lift of a quasicomplete, ISC surface of the form $(\Omega,\hat{e}_\infty,x_\infty)$, for some open subset $\Omega$ of $\Bbb{H}^2$.
\endproclaim
\proclabel{LabouriesCompactnessTheoremPGCII}
\remark[LabouriesCompactnessTheoremPGCII] Note, in particular, that the domain of the limit immersion need not be the whole of $\Bbb{H}^2$.
\newhead{Bibliography}[Bibliography]
\medskip
{\leftskip = 5ex \parindent = -5ex
\leavevmode\hbox to 4ex{\hfil \cite{Aron}}\hskip 1ex{Aronszajn N., A unique continuation theorem for elliptic differential equations or inequalities of the second order, {\sl J. Math. Pures Appl.}, {\bf 36}, (1957), 235--239}
\smallskip\leavevmode\hbox to 4ex{\hfil \cite{BarBegZegh}}\hskip 1ex{Barbot T., B\'eguin F., Zeghib A., Prescribing Gauss curvature of surfaces in $3$-dimensional spacetimes: application to the Minkowski problem in the Minkowski space, {\sl Ann. Inst. Fourier}, {\bf 61}, no. 2, (2011), 511–-591}
\smallskip\leavevmode\hbox to 4ex{\hfil \cite{BonMonSchI}}\hskip 1ex{Bonsante F., Mondello G., Schlenker J. M., A cyclic extension of the earthquake flow I, {\sl Geom. Topol.}, {\bf 17}, no. 1, (2013), 157--234}
\smallskip\leavevmode\hbox to 4ex{\hfil \cite{BonMonSchII}}\hskip 1ex{Bonsante F., Mondello G., Schlenker J. M., A cyclic extension of the earthquake flow II, {\sl Ann. Sci. Ec. Norm. Sup\'er.}, {\bf 48}, no. 4, (2015), 811–-859}
\smallskip\leavevmode\hbox to 4ex{\hfil \cite{Calabi}}\hskip 1ex{Calabi E., Improper affine hyperspheres of convex type and a generalisation of a theorem by K. J\"orgens, {\sl Mich. Math. J.}, {\bf 5}, (1958), 105--126}
\smallskip\leavevmode\hbox to 4ex{\hfil \cite{Girard}}\hskip 1ex{Girard P. R., The quaternion group and modern physics, {\sl Eur. J. Phys.}, {\bf 5}, (1984) 25--32}
\smallskip\leavevmode\hbox to 4ex{\hfil \cite{GromovPH}}\hskip 1ex{Gromov P. H., Pseudo-holomorphic curves in symplectic manifolds, {\sl Inv. Math.}, {\bf 82}, (1985), 307--347}
\smallskip\leavevmode\hbox to 4ex{\hfil \cite{Gromov}}\hskip 1ex{Gromov M., Foliated Plateau problem, Part II: Harmonic maps of foliations, {\sl Geom. Func. Anal.}, {\bf 1}, no. 3, (1991), 253--320}
\smallskip\leavevmode\hbox to 4ex{\hfil \cite{HarvLaws}}\hskip 1ex{Harvey R., Lawson H. B. Jr., Calibrated geometries, {\sl Acta. Math.}, {\bf 148}, (1982), 47--157}
\smallskip\leavevmode\hbox to 4ex{\hfil \cite{HarvLawsII}}\hskip 1ex{Harvey F. R., Lawson H. B., Pseudoconvexity for the special Lagrangian potential equation, {\sl Calc. Var. PDEs.}, {\bf 60}, (2021)}
\smallskip\leavevmode\hbox to 4ex{\hfil \cite{Jorgens}}\hskip 1ex{J\"orgens K., \"Uber die L\"osungen der Differentialgleichung $rt-s^2=1$, {\sl Math. Ann.}, {\bf 127}, (1954), 130--134}
\smallskip\leavevmode\hbox to 4ex{\hfil \cite{KokRossSajiUmeYam}}\hskip 1ex{Kokubu M., Rossman W., Saji K., Umehara M., Yamada K., Singularities of flat fronts in hyperbolic space, {\sl Pac. J. Math.}, {\bf 221}, no. 2, (2005), 303--351}
\smallskip\leavevmode\hbox to 4ex{\hfil \cite{KokRossUmeYam}}\hskip 1ex{Kokubu M., Rossman W., Umehara M., Yamada K., Flat fronts in hyperbolic space and their caustics, {\sl J. Soc. Math. Japan}, {\bf 59}, no. 1, (2007), 265--299}
\smallskip\leavevmode\hbox to 4ex{\hfil \cite{LabMP}}\hskip 1ex{Labourie F., M\'etriques prescrites sur le bord des vari\'et\'es hyperboliques de dimension $3$, {\sl J. Diff. Geom.}, {\bf 35}, no. 3, (1992), 609–-626}
\smallskip\leavevmode\hbox to 4ex{\hfil \cite{LabMA}}\hskip 1ex{Labourie F., Probl\`emes de Monge-Amp\`ere, courbes pseudo-holomorphes et laminations, {\sl Geom. Func. Anal.}, {\bf 7}, (1997), 496--534}
\smallskip\leavevmode\hbox to 4ex{\hfil \cite{LychRub}}\hskip 1ex{Lychagin V. V., Rubtsov V. N., Local classification of Monge-Amp\`ere differential equations (Russian), {\sl Dokl. Akad. Nauk SSSR}, {\bf 272}, no. 1, (1983), 34--38}
\smallskip\leavevmode\hbox to 4ex{\hfil \cite{LychRubChek}}\hskip 1ex{Lychagin V. V., Rubtsov V. N., Chekalov I. V., A classification of Monge-Amp\`ere equations, {\sl Ann. Sci. \'Ecole Norm. Sup.}, {\bf 26}, no. 3, (1993), 281--308}
\smallskip\leavevmode\hbox to 4ex{\hfil \cite{MartMil}}\hskip 1ex{Mart\'\i nez A., Mil\'an F., Flat fronts in hyperbolic $3$-space with prescribed singularities, {\sl Ann. Glob. Anal. Geom.}, {\bf 46}, (2014), 227--239}
\smallskip\leavevmode\hbox to 4ex{\hfil \cite{McDuffSal}}\hskip 1ex{McDuff D., Salamon D., {\sl J-holomorphic curves and quantum cohomology}, University Lecture Series, {\bf 6}, AMS, Providence, (1994)}
\smallskip\leavevmode\hbox to 4ex{\hfil \cite{Pogorelov}}\hskip 1ex{Pogorelov A. V., On the improper convex affine hyperspheres, {\sl Geom. Dedi.}, {\bf 1}, (1972), 33--46}
\smallskip\leavevmode\hbox to 4ex{\hfil \cite{Rub}}\hskip 1ex{Rubtsov V., Geometry of Monge-Ampère structures, in {\sl Nonlinear PDEs, their geometry, and applications}, Birkh\"auser/Springer, (2019), 95--156}
\smallskip\leavevmode\hbox to 4ex{\hfil \cite{Sch}}\hskip 1ex{Schlenker J. M., Hyperbolic manifolds with convex boundary, {\sl Inv. Math.}, {\bf 163}, (2006), 109--169}
\smallskip\leavevmode\hbox to 4ex{\hfil \cite{Spivak}}\hskip 1ex{Spivak M., {\sl A comprehensive introduction to differential geometry}, Vol. IV, Publish or Perish, (1999)}
\smallskip\leavevmode\hbox to 4ex{\hfil \cite{SmiAA}}\hskip 1ex{Smith G., An Arzela-Ascoli Theorem for Immersed Submanifolds, {\sl Ann. Fac. Sci. Toulouse Math.}, {\bf 16}, no. 4, (2007), 817--866}
\smallskip\leavevmode\hbox to 4ex{\hfil \cite{SmiSLC}}\hskip 1ex{Smith G., Special Lagrangian curvature, {\sl Math. Ann.}, {\bf 335}, no. 1, (2013), 57--95}
\smallskip\leavevmode\hbox to 4ex{\hfil \cite{SmiAS}}\hskip 1ex{Smith G., On the asymptotic geometry of finite-type k-surfaces in three-dimensional hyperbolic space, arXiv:1908.04834}
\smallskip\leavevmode\hbox to 4ex{\hfil \cite{TouLabWol}}\hskip 1ex{Toulisse J., Labourie F., Wolf M., Plateau Problems for Maximal Surfaces in Pseudo-Hyperbolic Spaces, to appear in {\sl Ann. Sci. \'Ec. Norm. Sup\'er}, arXiv:2006.12190}
\par}
%
%
%
%
\enddocument%
\ifloadreferences\else\closeout\references\fi%
\ifundefinedreferences\write16{There were undefined references.}\fi%
\ifchangedreferences\write16{References have changed.}\fi%
\end

%% file: references.tex
\global\def\_@citation@Aron{1}
\global\def\_@citation@BarBegZegh{2}
\global\def\_@citation@BonMonSchI{3}
\global\def\_@citation@BonMonSchII{4}
\global\def\_@citation@Calabi{5}
\global\def\_@citation@Girard{6}
\global\def\_@citation@GromovPH{7}
\global\def\_@citation@Gromov{8}
\global\def\_@citation@HarvLaws{9}
\global\def\_@citation@HarvLawsII{10}
\global\def\_@citation@Jorgens{11}
\global\def\_@citation@KokRossSajiUmeYam{12}
\global\def\_@citation@KokRossUmeYam{13}
\global\def\_@citation@LabMP{14}
\global\def\_@citation@LabMA{15}
\global\def\_@citation@LychRub{16}
\global\def\_@citation@LychRubChek{17}
\global\def\_@citation@MartMil{18}
\global\def\_@citation@McDuffSal{19}
\global\def\_@citation@Pogorelov{20}
\global\def\_@citation@Rub{21}
\global\def\_@citation@Sch{22}
\global\def\_@citation@Spivak{23}
\global\def\_@citation@SmiAA{24}
\global\def\_@citation@SmiSLC{25}
\global\def\_@citation@SmiAS{26}
\global\def\_@citation@TouLabWol{27}
\global\def\_@head@Introduction{1}
\global\def\_@subhead@Introduction{1.1}
\global\def\_@eqn@CurvaturePrescription{\relax \unhbox \voidb@x \hbox {{\relax \tenrm (1.1)}}}
\global\def\_@proc@IntroLabouriesCompactnessTheorem{1.1.1}
\global\def\_@rmk@IntroLabouriesCompactnessTheoremI{\relax \unhbox \voidb@x \hbox {1.1.1}}
\global\def\_@rmk@IntroLabouriesCompactnessTheoremII{\relax \unhbox \voidb@x \hbox {1.1.2}}
\global\def\_@proc@IntroLabouriesDichotomy{1.1.2}
\global\def\_@rmk@IntroLabouriesDichotomyI{\relax \unhbox \voidb@x \hbox {1.1.3}}
\global\def\_@subhead@Acknowledgements{1.2}
\global\def\_@head@QuaternionsAndBernsteinTypeTheorems{2}
\global\def\_@subhead@Quaternions{2.1}
\global\def\_@eqn@QuaternionRelations{\relax \unhbox \voidb@x \hbox {{\relax \tenrm (2.1)}}}
\global\def\_@eqn@GeneralQuaternionFormula{\relax \unhbox \voidb@x \hbox {{\relax \tenrm (2.2)}}}
\global\def\_@eqn@DefinitionOfQuaternionConjugate{\relax \unhbox \voidb@x \hbox {{\relax \tenrm (2.3)}}}
\global\def\_@eqn@ComplexConjugationsAnticommutes{\relax \unhbox \voidb@x \hbox {{\relax \tenrm (2.4)}}}
\global\def\_@eqn@RealityOfInnerProduct{\relax \unhbox \voidb@x \hbox {{\relax \tenrm (2.5)}}}
\global\def\_@eqn@QuaternionInnerProduct{\relax \unhbox \voidb@x \hbox {{\relax \tenrm (2.6)}}}
\global\def\_@eqn@FormulaForNormSquared{\relax \unhbox \voidb@x \hbox {{\relax \tenrm (2.7)}}}
\global\def\_@eqn@LengthIsMultiplicative{\relax \unhbox \voidb@x \hbox {{\relax \tenrm (2.8)}}}
\global\def\_@eqn@AntisymmetryI{\relax \unhbox \voidb@x \hbox {{\relax \tenrm (2.9)}}}
\global\def\_@eqn@VolumeForm{\relax \unhbox \voidb@x \hbox {{\relax \tenrm (2.10)}}}
\global\def\_@eqn@CoveringOfSO{\relax \unhbox \voidb@x \hbox {{\relax \tenrm (2.11)}}}
\global\def\_@proc@CoveringOfSO{2.1.1}
\global\def\_@eqn@KernelOfH{\relax \unhbox \voidb@x \hbox {{\relax \tenrm (2.12)}}}
\global\def\_@subhead@CompatibleComplexStructures{2.2}
\global\def\_@eqn@ComplexStructureCondition{\relax \unhbox \voidb@x \hbox {{\relax \tenrm (2.13)}}}
\global\def\_@eqn@JPreservesLength{\relax \unhbox \voidb@x \hbox {{\relax \tenrm (2.14)}}}
\global\def\_@proc@ComplexStructures{2.2.1}
\global\def\_@proc@LeftAndRightComplexStructures{2.2.2}
\global\def\_@rmk@LeftAndRightComplexStructures{\relax \unhbox \voidb@x \hbox {2.2.1}}
\global\def\_@proc@CharacterisationOfLeftComplexStructures{2.2.3}
\global\def\_@rmk@CharacterisationOfLeftComplexStructures{\relax \unhbox \voidb@x \hbox {2.2.2}}
\global\def\_@subhead@HomeomorphismsAndCompatibleQuaternionicStructures{2.3}
\global\def\_@eqn@CompatibilityCondition{\relax \unhbox \voidb@x \hbox {{\relax \tenrm (2.15)}}}
\global\def\_@eqn@HomomorphismsOfQuaternions{\relax \unhbox \voidb@x \hbox {{\relax \tenrm (2.16)}}}
\global\def\_@proc@HomomorphismsOfQuaternions{2.3.1}
\global\def\_@eqn@LeftQuaternionicStructure{\relax \unhbox \voidb@x \hbox {{\relax \tenrm (2.17)}}}
\global\def\_@eqn@RightQuaternionicStructure{\relax \unhbox \voidb@x \hbox {{\relax \tenrm (2.18)}}}
\global\def\_@proc@LeftAndRightQuaternionicStructures{2.3.2}
\global\def\_@rmk@LeftAndRightQuaternionicStructures{\relax \unhbox \voidb@x \hbox {2.3.1}}
\global\def\_@subhead@Calibrations{2.4}
\global\def\_@eqn@LeftComplexStructure{\relax \unhbox \voidb@x \hbox {{\relax \tenrm (2.19)}}}
\global\def\_@eqn@AssociatedSymplecticForm{\relax \unhbox \voidb@x \hbox {{\relax \tenrm (2.20)}}}
\global\def\_@eqn@Calibration{\relax \unhbox \voidb@x \hbox {{\relax \tenrm (2.21)}}}
\global\def\_@proc@CalibrationA{2.4.1}
\global\def\_@rmk@CalibrationA{\relax \unhbox \voidb@x \hbox {2.4.1}}
\global\def\_@proc@CalibrationB{2.4.2}
\global\def\_@subhead@TheMongeAmpereEquation{2.5}
\global\def\_@eqn@MongeAmpere{\relax \unhbox \voidb@x \hbox {{\relax \tenrm (2.22)}}}
\global\def\_@eqn@StandardComplexStructure{\relax \unhbox \voidb@x \hbox {{\relax \tenrm (2.23)}}}
\global\def\_@eqn@StandardMetric{\relax \unhbox \voidb@x \hbox {{\relax \tenrm (2.24)}}}
\global\def\_@eqn@QuaternionicStructure{\relax \unhbox \voidb@x \hbox {{\relax \tenrm (2.25)}}}
\global\def\_@eqn@BaseComplexStructure{\relax \unhbox \voidb@x \hbox {{\relax \tenrm (2.26)}}}
\global\def\_@eqn@SymplecticStructures{\relax \unhbox \voidb@x \hbox {{\relax \tenrm (2.27)}}}
\global\def\_@eqn@GraphPlane{\relax \unhbox \voidb@x \hbox {{\relax \tenrm (2.28)}}}
\global\def\_@eqn@TwoDimMatrixRelations{\relax \unhbox \voidb@x \hbox {{\relax \tenrm (2.29)}}}
\global\def\_@proc@FirstLexicon{2.5.1}
\global\def\_@proc@SolutionsOfMAAreJHolomorphicCurves{2.5.2}
\global\def\_@subhead@Positivity{2.6}
\global\def\_@eqn@VerticalSubspace{\relax \unhbox \voidb@x \hbox {{\relax \tenrm (2.30)}}}
\global\def\_@eqn@MinkowskiMetric{\relax \unhbox \voidb@x \hbox {{\relax \tenrm (2.31)}}}
\global\def\_@eqn@InvarianceOfM{\relax \unhbox \voidb@x \hbox {{\relax \tenrm (2.32)}}}
\global\def\_@proc@CharacterisationOfM{2.6.1}
\global\def\_@eqn@StabiliserSubgroup{\relax \unhbox \voidb@x \hbox {{\relax \tenrm (2.33)}}}
\global\def\_@proc@StabiliserSubgroup{2.6.2}
\global\def\_@rmk@StabiliserSubgroup{\relax \unhbox \voidb@x \hbox {2.6.1}}
\global\def\_@proc@SecondLexicon{2.6.3}
\global\def\_@subhead@PositivityIIAHolomorphicApproach{2.7}
\global\def\_@eqn@SpaceTimeCoordinatesI{\relax \unhbox \voidb@x \hbox {{\relax \tenrm (2.34)}}}
\global\def\_@eqn@SpaceTimeCoordinatesII{\relax \unhbox \voidb@x \hbox {{\relax \tenrm (2.35)}}}
\global\def\_@eqn@PushForwardOfMinkowskiMetric{\relax \unhbox \voidb@x \hbox {{\relax \tenrm (2.36)}}}
\global\def\_@eqn@Bilipschitz{\relax \unhbox \voidb@x \hbox {{\relax \tenrm (2.37)}}}
\global\def\_@proc@Bilipschitz{2.7.1}
\global\def\_@eqn@ParametrisationOfCP{\relax \unhbox \voidb@x \hbox {{\relax \tenrm (2.38)}}}
\global\def\_@proc@HolomorphicCharacterisationOfPositiveAndNullPlanes{2.7.2}
\global\def\_@subhead@BernsteinTypeTheorems{2.8}
\global\def\_@proc@JuergensQuaternionicSetting{2.8.1}
\global\def\_@proc@JuergensB{2.8.2}
\global\def\_@rmk@JurgensB{\relax \unhbox \voidb@x \hbox {2.8.1}}
\global\def\_@proc@VolkovVladimirovaSasaki{2.8.3}
\global\def\_@rmk@VolkovVladimirovaSasaki{\relax \unhbox \voidb@x \hbox {2.8.2}}
\global\def\_@subhead@BernsteinTypeTheoremsII{2.9}
\global\def\_@eqn@BoundaryCondition{\relax \unhbox \voidb@x \hbox {{\relax \tenrm (2.39)}}}
\global\def\_@proc@BoundaryJuergensQuaternionicSetting{2.9.1}
\global\def\_@proc@BoundaryJuergens{2.9.2}
\global\def\_@subhead@NonQuadraticSolutions{2.10}
\global\def\_@eqn@FatalCounterExample{\relax \unhbox \voidb@x \hbox {{\relax \tenrm (2.40)}}}
\global\def\_@eqn@Ansatz{\relax \unhbox \voidb@x \hbox {{\relax \tenrm (2.41)}}}
\global\def\_@eqn@VerticalSpaceInSpaceTime{\relax \unhbox \voidb@x \hbox {{\relax \tenrm (2.42)}}}
\global\def\_@eqn@CounterExampleBoundaryConditionA{\relax \unhbox \voidb@x \hbox {{\relax \tenrm (2.43)}}}
\global\def\_@eqn@CounterExampleBoundaryConditionB{\relax \unhbox \voidb@x \hbox {{\relax \tenrm (2.44)}}}
\global\def\_@eqn@DoublingProperty{\relax \unhbox \voidb@x \hbox {{\relax \tenrm (2.45)}}}
\global\def\_@eqn@FirstComponentOfCounterExample{\relax \unhbox \voidb@x \hbox {{\relax \tenrm (2.46)}}}
\global\def\_@eqn@FormulaForF{\relax \unhbox \voidb@x \hbox {{\relax \tenrm (2.47)}}}
\global\def\_@eqn@TransformedFormulaForF{\relax \unhbox \voidb@x \hbox {{\relax \tenrm (2.48)}}}
\global\def\_@eqn@OneForm{\relax \unhbox \voidb@x \hbox {{\relax \tenrm (2.49)}}}
\global\def\_@head@MongeAmpereStructuresAndKSurfaces{3}
\global\def\_@subhead@MongeAmpereStructures{3.1}
\global\def\_@rmk@AlternativePerspective{\relax \unhbox \voidb@x \hbox {3.1.1}}
\global\def\_@eqn@JHolomorphicCurve{\relax \unhbox \voidb@x \hbox {{\relax \tenrm (3.1)}}}
\global\def\_@eqn@CurtainSurfaceTheorem{\relax \unhbox \voidb@x \hbox {{\relax \tenrm (3.2)}}}
\global\def\_@proc@CurtainSurfaceTheorem{3.1.1}
\global\def\_@proc@CurtainSurfaceTheoremI{3.1.2}
\global\def\_@proc@CurtainSurfaceTheoremII{3.1.3}
\global\def\_@eqn@PolynomialApproximation{\relax \unhbox \voidb@x \hbox {{\relax \tenrm (3.3)}}}
\global\def\_@eqn@QVanishes{\relax \unhbox \voidb@x \hbox {{\relax \tenrm (3.4)}}}
\global\def\_@eqn@FormulaForP{\relax \unhbox \voidb@x \hbox {{\relax \tenrm (3.5)}}}
\global\def\_@subhead@Compactness{3.2}
\global\def\_@proc@Precompactness{3.2.1}
\global\def\_@rmk@Precompactness{\relax \unhbox \voidb@x \hbox {3.2.1}}
\global\def\_@proc@QuasiMaximumLemma{3.2.2}
\global\def\_@proc@Compactness{3.2.3}
\global\def\_@proc@LabouriesCompactnessTheorem{3.2.4}
\global\def\_@subhead@ApplicationsIKSurfaces{3.3}
\global\def\_@eqn@HorizontalAndVerticalSubbundles{\relax \unhbox \voidb@x \hbox {{\relax \tenrm (3.6)}}}
\global\def\_@eqn@DecompositionOfTangentBundle{\relax \unhbox \voidb@x \hbox {{\relax \tenrm (3.7)}}}
\global\def\_@eqn@MAWAndV{\relax \unhbox \voidb@x \hbox {{\relax \tenrm (3.8)}}}
\global\def\_@eqn@DefinitionOfWedgeProduct{\relax \unhbox \voidb@x \hbox {{\relax \tenrm (3.9)}}}
\global\def\_@eqn@ComplexStructureOverV{\relax \unhbox \voidb@x \hbox {{\relax \tenrm (3.10)}}}
\global\def\_@eqn@MAComplexStructure{\relax \unhbox \voidb@x \hbox {{\relax \tenrm (3.11)}}}
\global\def\_@eqn@MApplicationI{\relax \unhbox \voidb@x \hbox {{\relax \tenrm (3.12)}}}
\global\def\_@proc@IntegrabilityI{3.3.1}
\global\def\_@eqn@TangentSpaceOfNGamma{\relax \unhbox \voidb@x \hbox {{\relax \tenrm (3.13)}}}
\global\def\_@eqn@PrescribedCurvature{\relax \unhbox \voidb@x \hbox {{\relax \tenrm (3.14)}}}
\global\def\_@proc@DescriptionOfMongeAmpereSurfacesI{3.3.2}
\global\def\_@proc@LabouriesCompactnessTheoremPGCI{3.3.3}
\global\def\_@rmk@LabouriesCompactnessTheoremPGCI{\relax \unhbox \voidb@x \hbox {3.3.1}}
\global\def\_@subhead@ApplicationsIIIsometricImmersions{3.4}
\global\def\_@eqn@TrivialisationI{\relax \unhbox \voidb@x \hbox {{\relax \tenrm (3.15)}}}
\global\def\_@eqn@TrivialisationII{\relax \unhbox \voidb@x \hbox {{\relax \tenrm (3.16)}}}
\global\def\_@eqn@HorizontalAndVerticalSubbundlesII{\relax \unhbox \voidb@x \hbox {{\relax \tenrm (3.17)}}}
\global\def\_@eqn@DecompositionOfTIsom{\relax \unhbox \voidb@x \hbox {{\relax \tenrm (3.18)}}}
\global\def\_@eqn@DefinitionOfWAndVII{\relax \unhbox \voidb@x \hbox {{\relax \tenrm (3.19)}}}
\global\def\_@eqn@DefinitionOfJII{\relax \unhbox \voidb@x \hbox {{\relax \tenrm (3.20)}}}
\global\def\_@eqn@DefinitionMII{\relax \unhbox \voidb@x \hbox {{\relax \tenrm (3.21)}}}
\global\def\_@eqn@PropertiesOfAlpha{\relax \unhbox \voidb@x \hbox {{\relax \tenrm (3.22)}}}
\global\def\_@eqn@DefinitionOfCurtainSurfaceII{\relax \unhbox \voidb@x \hbox {{\relax \tenrm (3.23)}}}
\global\def\_@eqn@DerivativeOfCurtainII{\relax \unhbox \voidb@x \hbox {{\relax \tenrm (3.24)}}}
\global\def\_@eqn@GaussLiftII{\relax \unhbox \voidb@x \hbox {{\relax \tenrm (3.25)}}}
\global\def\_@eqn@PrescribedCurvatureII{\relax \unhbox \voidb@x \hbox {{\relax \tenrm (3.26)}}}
\global\def\_@proc@DescriptionOfMongeAmpereSurfacesII{3.4.2}
\global\def\_@eqn@DefinitionOfPhiFinalExample{\relax \unhbox \voidb@x \hbox {{\relax \tenrm (3.27)}}}
\global\def\_@eqn@DefinitionOfHatEm{\relax \unhbox \voidb@x \hbox {{\relax \tenrm (3.28)}}}
\global\def\_@proc@LabouriesCompactnessTheoremPGCII{3.4.3}
\global\def\_@rmk@LabouriesCompactnessTheoremPGCII{\relax \unhbox \voidb@x \hbox {3.4.1}}
\global\def\_@head@Bibliography{4}